\documentclass[11pt]{article}
\usepackage[english]{babel}
\usepackage{amsmath,amsthm,amssymb}
\usepackage{mathrsfs}
\usepackage{bbm}
\usepackage{amsfonts}
\usepackage[margin=0.8in]{geometry}
\usepackage{graphicx}
\usepackage{caption}
\usepackage{subcaption}
\usepackage{wasysym}
\usepackage{enumerate}
\usepackage{algorithm2e}
\usepackage{tabularx}
\usepackage{color}
\usepackage{wrapfig}
\usepackage{todonotes}

\newtheorem{thm}{Theorem}[section]
\newtheorem{ass}[thm]{Assumption}
\newtheorem{cor}[thm]{Corollary}
\newtheorem{lem}[thm]{Lemma}
\newtheorem{prop}[thm]{Proposition}
\newtheorem*{hyp*}{Hypothesis}
\theoremstyle{definition}
\newtheorem{defn}[thm]{Definition}

\theoremstyle{rem}
\newtheorem{rem}[thm]{Remark}

\numberwithin{equation}{section}

\newcommand{\R}{\mathbb R}

\newcommand{\bbF}{\mathbb F}

\newcommand{\mcA}{\mathcal{A}}

\newcommand{\mcL}{\mathcal L}

\newcommand{\mcF}{\mathcal F}

\newcommand{\mcT}{\mathcal T}

\newcommand{\mcM}{\mathcal M}
\newcommand{\mcH}{\mathcal H}
\newcommand{\mcK}{\mathcal K}
\newcommand{\mcU}{\mathcal U}
\newcommand{\mcS}{\mathcal S}
\newcommand{\mcY}{\mathcal Y}
\newcommand{\mcZ}{\mathcal Z}

\newcommand{\bigS}{\mathfrak S}
\newcommand{\bigH}{\mathfrak H}

\newcommand{\E}{\mathbb{E}}
\newcommand{\Prob}{\mathbb{P}}

\DeclareMathOperator*{\esssup}{ess\,sup}

\newcommand{\ett}{\mathbbm{1}}

\newcommand{\cadlag}{c\`adl\`ag~}
\newcommand{\caglad}{c\`agl\`ad~}

\newcommand{\trace}{{\rm Tr}}

\newcommand{\ie}{\textit{i.e.\ }}
\newcommand{\eg}{\textit{e.g.\ }}
\newcommand{\etal}{\textit{et.~al.\ }}

\begin{document}

\title{Probabilistic Representation of Viscosity Solutions to Quasi-Variational Inequalities with Non-Local Drivers\footnote{This work was supported by the Swedish Energy Agency through grant number 48405-1}}

\author{Magnus Perninge\footnote{M.\ Perninge is with the Department of Physics and Electrical Engineering, Linnaeus University, V\"axj\"o,
Sweden. e-mail: magnus.perninge@lnu.se.}} %
\maketitle
\begin{abstract}
We consider quasi-variational inequalities (QVIs) with general non-local drivers and related systems of reflected backward stochastic differential equations (BSDEs) in a Brownian filtration. We show existence and uniqueness of viscosity solutions to the QVIs by first considering the standard (local) setting and then applying a contraction argument. In addition, the contraction argument yields existence and uniqueness of solutions to the related systems of reflected BSDEs and extends the theory of probabilistic representations of PDEs in terms of BSDEs to our specific setting.
\end{abstract}

\section{Introduction}
We consider existence and uniqueness of viscosity solutions to the quasi-variational inequality (QVI)
\begin{align}\label{ekv:var-ineq}
\begin{cases}
  \min\{v(t,x)-\mcM v(t,x),-v_t(t,x)-\mcL v(t,x)-f(t,x,v(t,\cdot),\sigma^\top(t,x)\nabla_x v(t,x))\}=0,\\
  \quad\forall (t,x)\in[0,T)\times \R^d \\
  v(T,x)=\psi(x),
\end{cases}
\end{align}
where for each $(t,x,z)\in [0,T]\times\R^n\times\R^d$, the map $g\mapsto f(t,x,g,z):C(\R^n\to\R)\to\R$ is a functional, $\mcM v(t,x):=\sup_{b\in U}\{v(t,\Gamma(t,x,b))-\ell(t,x,b)\}$ and
\begin{align}
  \mcL:=\sum_{j=1}^d a_j(t,x)\frac{\partial}{\partial x_j}+\frac{1}{2}\sum_{i,j=1}^d (\sigma\sigma^\top(t,x))_{i,j}\frac{\partial^2}{\partial x_i\partial x_j}
\end{align}
is the infinitesimal generator related to the SDE
\begin{align*}
X^{t,x}_s=x + \int_t^s a(r,X^{t,x}_r)dr + \int_t^s \sigma(r,X^{t,x}_r)dW_r.
\end{align*}
We take a probabilistic approach and extend the Feynman-Kac framework by relating solutions to \eqref{ekv:var-ineq} to those of the system of reflected BSDEs
\begin{align}\label{ekv:syst-bsde-new}
\begin{cases}
  Y^{t,x}_s=\psi(X^{t,x}_T)+\int_s^T f(r,X^{t,x}_r,\bar Y(r,\cdot),Z^{t,x}_r)dr-\int_s^T Z^{t,x}_r dW_r+ K^{t,x}_T-K^{t,x}_s,\quad\forall s\in[t,T], \\
  Y^{t,x}_s\geq\mcM\bar Y(s,X^{t,x}_s),\:\forall s\in[t,T] \quad{\rm and}\quad \int_t^T(Y^{t,x}_s-\mcM \bar Y(s,X^{t,x}_s))dK^{t,x}_s=0,
\end{cases}
\end{align}
where $\bar Y$ is a continuous, deterministic function such that $\bar Y(t,x)=Y^{t,x}_t$, $\Prob$-a.s.

It is well known that, under suitable conditions on the involved parameters, value functions to impulse control problems are solutions (in viscosity sense) to standard QVIs when the driving noise process is a Brownian motion (see the seminal work in~\cite{BensLionsImpulse}) and to so called quasi-integrovariational inequalities when the driving noise is a general L\'evy process~\cite{OksenSulemBok}.

Our primary motivation for studying the QVIs of the more general type in \eqref{ekv:var-ineq} and the corresponding systems of reflected BSDEs \eqref{ekv:syst-bsde-new} is this close connection to stochastic impulse control that we exploit in the case of a local driver $f:[0,T]\times\R^n\times\R\times\R^d\to\R$ as a means to derive an intermediate result. Moreover, existence and uniqueness of solutions to \eqref{ekv:syst-bsde-new} appears to be a cornerstone in attempts to prove existence of solutions to systems of doubly reflected BSDEs through penalization (see \cite{BollanSWG1,DjehicheSWG2} for a complete result in a less general setup). Since the unique solution to the system of doubly reflected BSDEs in \cite{DjehicheSWG2} appear as the value function to a switching game\footnote{An optimal switching problem is a type of impulse control problem where the operator of a system seeks to maximize the revenue by switching between a set of operation modes (see~\cite{CarmLud} for a complete definition).} in \cite{HamMM_SWG}, it is natural to assume that our results can lead to an extension of the theory of stochastic differential games of impulse control (see \eg \cite{Cosso13}). Given the generality of our framework and in particulary the type of non-locality in our driver it is also plausible that our results opens up new perspectives for dealing with more general impulse control problems, for instance mean-field type control problems involving impulses.

An alternative Feynman-Kac representation for solutions to standard QVIs was proposed in \cite{KharroubiQVI}, where the solution to a QVI is related to the minimal solution of a BSDE driven by a Brownian motion and a Poisson random measure with a constraint on the jump component. It should be noted that a similar approach was take in \cite{KharroubiPham15} to derive a probabilistic representation for solutions to nonlinear integro-partial differential equations (IPDEs) of Hamilton-Jacobi-Bellman type.

Systems of reflected BSDEs related to non-Markovian optimal switching problems were independently considered in~\cite{HamZhang} and \cite{HuTang}. Based on a contraction argument, the authors of \cite{ChassElieKharr} later extended the existence and uniqueness results for such systems of reflected BSDEs, removing certain monotonicity assumptions on the drivers. In the Markovian framework, the corresponding Feynman-Kac representation was derived in \cite{Morlais13} under a monotonicity assumption on the driver and extended to the general setting in \cite{HamMnif19}.

Despite the surge of interest in switching problems and the related systems of reflected BSDEs, for a long while no work appeared that connected impulse control problems to reflected BSDEs. However, recently \cite{P_ImpBSDE2021} considered a sequential system of reflected BSDEs related to non-Markovian impulse control. The intermediate result in the present work extends that of \cite{P_ImpBSDE2021} in the Markovian framework and bridges the gap left in terms of a probabilistic representation of the corresponding QVIs.

As noted above, our intermediate step assumes that $f$ is a local operator and seeks a solution the system of reflected BSDEs
\begin{align}\label{ekv:seq-bsde-loc}
\begin{cases}
  Y^{t,x}_s=\psi(X^{t,x}_T)+\int_s^T f(r,X^{t,x}_r,Y^{t,x}_r,Z^{t,x}_r)dr-\int_s^T Z^{t,x}_r dW_r+ K^{t,x}_T-K^{t,x}_s,\quad\forall s\in[t,T], \\
  Y^{t,x}_s\geq\mcM\bar Y(s,X^{t,x}_s),\:\forall s\in[t,T] \quad{\rm and}\quad \int_t^T(Y^{t,x}_s-\mcM \bar Y(s,X^{t,x}_s))dK^{t,x}_s=0.
\end{cases}
\end{align}
Existence is achieved through a Piccard iteration approach and uniqueness is shown by a verification argument, relating a solution to \eqref{ekv:seq-bsde-loc} to the value function of a stochastic impulse control problem. Moreover, we show that if $(Y,Z,K)$ solves \eqref{ekv:seq-bsde-loc}, letting $v:[0,T]\times\R^n\to\R$ be such that $v(t,x):=Y^{t,x}_t$, $\Prob$-a.s., then $v$ is a viscosity solution to the QVI
\begin{align}\label{ekv:var-ineq-loc}
\begin{cases}
  \min\{v(t,x)-\mcM v(t,x),-v_t(t,x)-\mcL v(t,x)-f(t,x,v(t,x),\sigma^\top(t,x)\nabla_x v(t,x))\}=0,\\
  \quad\forall (t,x)\in[0,T)\times \R^d \\
  v(T,x)=\psi(x),
\end{cases}
\end{align}
and (in an Appendix) we show uniqueness of solutions to \eqref{ekv:var-ineq-loc}. To extend these results to the non-local setting we consider a sequence of interconnected quasi-variational inequalities
\begin{align}\label{ekv:var-ineq-k}
\begin{cases}
  \min\{v^k(t,x)-\mcM v^{k}(t,x),-v^k_t(t,x)-\mcL v^k(t,x)-f(t,x,v^{k-1}(t,\cdot),\sigma(t,x)\nabla_x v^k(t,x))\}=0,\\
  \quad\forall (t,x)\in[0,T)\times \R^d, \\
  v^k(T,x)=\psi(x),
\end{cases}
\end{align}
with $v_0 \equiv 0$ and find a norm under which the map $\Phi$ mapping $v_k$ to $v_{k+1}$ is a contraction.

The remained of the article is organised as follows. In the next section we set the notation and state the assumptions that hold throughout. In addition, we give some preliminary results that are repeatedly referred in the article. Then, in Section~\ref{sec:impulse-control} we turn to the local setting before we, in the following section, derive the complete result. Uniqueness of solutions to~\eqref{ekv:var-ineq-loc} appears rudimentary and resembles the corresponding results in \cite{Morlais13}. However, since our setting is fundamentally different and for the sake of completeness, a uniqueness proof through viscosity comparison is included as an appendix.

\section{Preliminaries\label{sec:prel}}
\subsection{Notation}
We let $(\Omega,\mcF,\Prob)$ be a complete probability space on which lives a $d$-dimensional Brownian motion $W$. We denote by $\bbF:=(\mcF_t)_{0\leq t\leq T}$ the augmented natural filtration of $W$ and for $t\in[0,T]$ we let $\bbF^t:=(\mcF^t_s)_{t\leq s\leq T}$ denote the augmented natural filtration generated by $(W_s-W_t:t\leq s\leq T)$.\\

\noindent Throughout, we will use the following notation, where $n\geq 1$ is the dimension of the state-space:
\begin{itemize}
\item We let $\Pi^g$ denote the set of all functions $\varphi:[0,T]\times\R^n\to\R$ that are of polynomial growth in $x$, \ie there are constants $C,\rho>0$ such that $|\varphi(t,x)|\leq C(1+|x|^\rho)$ for all $(t,x)\in [0,T]\times\R^n$, and let $\Pi^g_c$ be the subset of jointly continuous functions.
  \item For $p\geq 1$ and $t\in [0,T]$, we let $\mcS^{p}_{\text{cag},t}$ be the set of all $\R$-valued, $\bbF^t$-progressively measurable \caglad processes $(Z_s: s\in [t,T])$ such that $\|Z\|_{\mcS^p_t}:=\E\big[\sup_{s\in[t,T]} |Z_s|^p\big]<\infty$. Moreover, we let $\mcS^p_{t}$ be the subset of continuous processes and $\mcS^p_{t,i}$ be the subset of continuous and non-decreasing processes with $Z_t=0$.
  \item We let $\mcH^{p}_t$ denote the set of all $\R^d$-valued $\bbF^t$-progressively measurable processes $(Z_s: s\in[t,T])$ such that $\|Z\|_{\mcH^p_t}:=\E\big[\big(\int_t^T |Z_s|^2 ds\big)^{p/2}\big]^{1/p}<\infty$.
  \item We let $\bigS^{p}$ be the set of all maps $Z:\cup_{t\in[0,T]}[t,T]\times \Omega\times \{t\}\times \R^n\to \R:(s,\omega,t,x)\mapsto Z^{t,x}_s(\omega)$ such that $Z^{t,x}\in \mcS^p_t$ and there is a $v\in\Pi^g_c$ such that $v(t,x)=Z^{t,x}_t$, $\Prob$-a.s., for all $(t,x)\in [0,T]\times\R^n$. Moreover, given $Z\in\bigS^p$, we let $\bar Z$ denote this deterministic function, so that $\bar Z(t,x)=v(t,x)$ for all $(t,x)\in[0,T]\times \R^n$.
  \item We let $\bigS^p_i$ be the subset of $\bigS^{p}$ with all maps $Z$ such that $Z^{t,x}\in \mcS^p_{t,i}$.
  \item We let $\bigH^{p}$ denote the set of all maps $Z:\cup_{t\in[0,T]}[t,T]\times \Omega\times \{t\}\times \R^n\to \R:(s,\omega,t,x)\mapsto Z^{t,x}_s(\omega)$ such that $Z^{t,x}\in \mcH^p_t$.
  \item We let $\mcT$ be the set of all $\bbF$-stopping times and for each $\eta\in\mcT$ we let $\mcT_\eta$ be the subset of stopping times $\tau$ such that $\tau\geq \eta$, $\Prob$-a.s.
  \item For $t\in[0,T]$, we let $\mcU_t$ be the set of all $u=(\tau_j,\beta_j)_{1\leq j\leq N}$, where $(\tau_j)_{j=1}^\infty$ is a non-decreasing sequence of $\bbF^t$-stopping times in $\mcT_t$, $\beta_j$ is a $\mcF^t_{\tau_j}$-measurable r.v.~taking values in the compact set $U$ and $N:=\max\{j:\tau_j<T\}$, such that $\Xi^{t,x;u}_T\in L^2(\Omega,\mcF^t_T,\Prob)$, where $\Xi^{t,x;u}_T$ is the total cost of impulses (see \eqref{ekv:Xi-def}).
  \item For $u\in\mcU$, we let $[u]_{j}:=(\tau_i,\beta_i)_{1\leq i\leq N\wedge j}$.
  \item For $t\in[0,T]$, we let $\mcA_t$ denote the set of all $[-1,1]^d$-valued, $\bbF^t$-progressively measurable processes $(\alpha_s:t\leq s\leq T)$ and set $\mcA:=\mcA_0$.
  \item For $t\in[0,T]$, we define the composition $\oplus_t$ of $\alpha^1\in\mcA$ and $\alpha^2\in\mcA_t$ as $(\alpha^1\oplus_t\alpha^2)_s:=\ett_{[0,t)}(s)\alpha^1_s+\ett_{[t,T]}(s)\alpha^2_s$.
\end{itemize}

We also mention that, unless otherwise specified, all inequalities between random variables are to be interpreted in the $\Prob$-a.s.~sense.

\subsection{Assumptions}

Throughout, we make the following assumptions on the parameters in the reward functional where $\rho>0$ is a fixed constant:
\begin{ass}\label{ass:on-coeff}
\begin{enumerate}[i)]
  \item\label{ass:on-coeff-f} We assume that $f:[0,T]\times \R^n\times C(\R^n\to\R)\times\R^{d}\to \R$ is such that for $v\in\Pi^g_c$, the map $(t,x)\mapsto f(t,x,v(t,\cdot),z)$ is jointly continuous, uniformly in $z$, $f$ is of polynomial growth in $x$, \ie there is a $C_f>0$ such that
  \begin{align*}
    |f(t,x,0,0)|\leq C_f(1+|x|^\rho)
  \end{align*}
  and that there are constants $k_f,K_\Gamma>0$ such that for any $t\in[0,T]$, $x\in\R^n$, $g,\tilde g\in C(\R^n\to\R)$ and $z,\tilde z\in\R^{d}$ we have
  \begin{align*}
    |f(t,x,\tilde g,\tilde z)-f(t,x,g,z)|&\leq k_f(\sup_{x'\in \Lambda_f(|x|)}|\tilde g(x')-g(x)|+|\tilde z-z|),
  \end{align*}
  where for each $\gamma\in\R_+$, $\Lambda_f(\gamma):=\{x\in\R^n:\|x\|\leq \gamma\vee K_\Gamma\}$ is the closed ball of radius $\gamma\vee K_\Gamma$ centered at the origin.
  \item\label{ass:on-coeff-psi} The terminal reward $\psi:\R^n\to\R$ is continuous and satisfies the growth condition
  \begin{align*}
    |\psi(x)|\leq C_\psi(1+|x|^\rho)
  \end{align*}
  for some $C_\psi>0$.
  \item\label{ass:on-coeff-ell} The intervention cost $\ell:[0,T]\times \R^n\times U\to \R_+$ is jointly continuous, of polynomial growth and bounded from below, \ie
  \begin{align*}
    \ell(t,x,b)\geq\delta >0,
  \end{align*}
  \item\label{ass:on-coeff-@end} For each $(x,b)\in\R^n\times U$ we have
  \begin{align*}
    \psi(x)>\psi(\Gamma(T,x,b))-\ell(t,x,b).
  \end{align*}
\end{enumerate}
\end{ass}

Moreover, we make the following assumptions on the coefficients of the forward SDE (and its impulsively controlled counterpart in Section~\ref{sec:impulse-control}):

\begin{ass}\label{ass:onSFDE}
For any $t,t'\in [0,T]$, $b\in U$ and $x,x'\in\R^n$ we have:
\begin{enumerate}[i)]
  \item\label{ass:onSFDE-Gamma} The function $\Gamma:[0,T]\times\R^n\times U\to\R^n$ is jointly continuous and satisfies the growth condition
  \begin{align}\label{ekv:imp-bound}
    |\Gamma(t,x,b)|\leq K_\Gamma\vee |x|.
  \end{align}
  \item\label{ass:onSFDE-a-sigma} The coefficients $a:[0,T]\times\R^n\to\R^{n}$ and $\sigma:[0,T]\times\R^n\to\R^{n\times d}$ are jointly continuous and satisfy the growth conditions
  \begin{align*}
    |a(t,x)|+|\sigma(t,x)|&\leq C_{a,\sigma}(1+|x|),
  \end{align*}
  for some $C_{a,\sigma}>0$ and the Lipschitz continuity
  \begin{align*}
    |a(t,x)-a(t,x')|+|\sigma(t,x)-\sigma(t,x')|&\leq k_{a,\sigma}|x'-x|,
  \end{align*}
  for some $k_{a,\sigma}>0$.
\end{enumerate}
\end{ass}

\subsection{Viscosity solutions}
We define the upper, $v^*$, and lower, $v_*$ semi-continuous envelope of a function $v$ as
\begin{align*}
v^*(t,x):=\limsup_{(t',x')\to(t,x),\,t'<T}v(t',x')\quad {\rm and}\quad v_*(t,x):=\liminf_{(t',x')\to(t,x),\,t'<T}v(t',x')
\end{align*}
Next we introduce the notion of a viscosity solution using the limiting parabolic superjet $\bar J^+v$ and subjet $\bar J^-v$ of a function $v$ (see pp. 9-10 of \cite{UsersGuide} for a definition):
\begin{defn}\label{def:visc-sol-jets}
Let $v$ be a locally bounded function from $[0,T]\times \R^n$ to $\R$. Then,
\begin{enumerate}[a)]
  \item It is referred to as a viscosity supersolution (resp. subsolution) to \eqref{ekv:var-ineq} if it is l.s.c.~(resp u.s.c.) and satisfies:
  \begin{enumerate}[i)]
    \item $v(T,x)\geq \psi(x)$ (resp. $v(T,x)\leq \psi(x)$)
    \item For any $(t,x)\in [0,T)\times\R^d$ and $(p,q,X)\in \bar J^- v(t,x)$ (resp. $\bar J^+ v(t,x)$) we have
    \begin{align*}
      \min\Big\{&v(t,x)-\mcM v(t,x),-p-H(t,x,v(t,x),q,X,a)\Big\}\geq 0
    \end{align*}
    (resp.
    \begin{align*}
      \min\Big\{&v(t,x)-\mcM v(t,x),-p-H(t,x,v(t,x),q,X,a)\Big\}\leq 0).
    \end{align*}
  \end{enumerate}
  \item It is called a viscosity solution to \eqref{ekv:var-ineq} if $v_*$ is a supersolution and $v^*$ is a subsolution.
\end{enumerate}
\end{defn}

We will sometimes use the following alternative definition of viscosity supersolutions (resp. subsolutions):
\begin{defn}\label{def:visc-sol-dom}
  A l.s.c.~(resp. u.s.c.) function $v$ is a viscosity supersolution (subsolution) to \eqref{ekv:var-ineq} if $v(T,x)\leq \psi(x)$ (resp. $\geq \psi(x)$) and whenever $\varphi\in C([0,T]\times\R^d\to\R)$ is such that $\varphi(t,x)=v(t,x)$ and $\varphi-v$ has a local maximum (resp. minimum) at $(t,x)$, then
  \begin{align*}
    \min\big\{&v(t,x)-\mcM v(t,x),-\varphi_t(t,x)-H(t,x,v(t,x),D\varphi(t,x),D^2\varphi(t,x),a)\big\}\geq 0\:(\leq 0).
  \end{align*}
\end{defn}


\subsection{Reflected BSDEs and obstacle problems}
We will make extensive use of the following classical result:

\begin{thm}\label{thm:ElKaroui}(El Karoui \etal \cite{ElKaroui1})
Let $h:[0,T]\times\R^n\to\R$ and $\psi:\R^n\to\R$ be (jointly) continuous and of polynomial growth with $h(T,x)\leq\Psi(x)$ for all $x\in\R^n$ and assume that $f:[0,T]\times\R^n\times\R\times\R^d\to\R$ ($(t,x,y,z)\mapsto f(t,x,y,z)$) is jointly continuous in $(t,x)$ uniformly in $(y,z)$, Lipschitz continuous in $(z,y)$ uniformly in $(t,x)$ and such that $f(\cdot,\cdot,0,0)\in\Pi^g$. Then, for each $(t,x)\in[0,T]\times\R^n$, there is a unique triple $(Y^{t,x},Z^{t,x},K^{t,x})\in \mcS^2_{t}\times\mcH^2_t\times\mcS^2_{t,i}$ such that
\begin{align}\label{ekv:rbsde-ElKaroui}
\begin{cases}
  Y^{t,x}_s=\psi(X^{t,x}_T)+\int_s^T f(r,X^{t,x}_r,Y^{t,x}_r,Z^{t,x}_r)dr-\int_s^T Z^{t,x}_r dW_r+ K^{t,x}_T-K^{t,x}_s,\quad\forall s\in[t,T], \\
  Y^{t,x}_s\geq h(s,X^{t,x}_s),\:\forall s\in[t,T] \quad{\rm and}\quad\int_t^T(Y^{t,x}_s-h(s,X^{t,x}_s))dK^{t,x}_s=0.
\end{cases}
\end{align}
Moreover, there is a $v\in \Pi^g_c$ such that $Y^{t,x}_s=v(s,X^{t,x}_s)$, $\Prob$-a.s.~for all $s\in[t,T]$, and $v$ is the unique viscosity solution in $\Pi^g$ to the following obstacle problem:
\begin{align}\label{ekv:obst-prob-ElKaroui}
\begin{cases}
  \min\{v(t,x)-h(t,x),-v_t(t,x)-\mcL v(t,x)-f(t,x,v(t,x),\sigma^\top(t,x)\nabla_x v(t,x))\}=0,\\
  \quad\forall (t,x)\in[0,T)\times \R^d \\
  v(T,x)=\psi(x).
\end{cases}
\end{align}
Furthermore\footnote{Throughout, $C$ will denote a generic positive constant that may change value from line to line.},
\begin{align}\label{ekv:ElKaroui-bound}
\|Y^{t,x}\|_{\mcS^2_t}^2+\|Z^{t,x}\|_{\mcH^2_t}^2+\|K^{t,x}\|_{\mcS^2_t}^2&\leq C\E\Big[|\psi(X^{t,x}_T)|^{2}+\int_t^T|f(s,X^{t,x}_s,0,0)|^2ds+\sup_{s\in[t,T]}|(h(s,X^{t,x}_s)^+|^2\Big].
\end{align}
In addition, $Y$ can be interpreted as the Snell envelope in the following way
\begin{equation*}
      Y^{t,x}_s=\esssup_{\tau\in\mcT_s}\E\bigg[\int_s^\tau f(r,Y^{t,x}_r,Z^{t,x}_r)dr+h(\tau,X^{t,x}_\tau)\ett_{[\tau<T]}+\psi(X^{t,x}_T) \ett_{[\tau=T]}\Big|\mcF^t_s\bigg]
\end{equation*}
and with $D^{t,x}_s:=\inf\{r\geq s: Y^{t,x}_r=h(r,X^{t,x}_r)\}\wedge T$ we have the representation
\begin{equation*}
      Y_s=\E\bigg[\int_s^{D^{t,x}_s} f(r,X^{t,x}_r,Y^{t,x}_r,Z^{t,x}_r)dr+h({D^{t,x}_s},X^{t,x}_{{D^{t,x}_s}})\ett_{[{D^{t,x}_s}<T]}+\psi(X^{t,x}_T) \ett_{[{D^{t,x}_s}=T]}\Big|\mcF^t_s\bigg]
\end{equation*}
and $K^{t,x}_{D^{t,x}_s}-K^{t,x}_s=0$, $\Prob$-a.s.

Finally, if $(\tilde Y,\tilde Z,\tilde K)$ is the solution to the reflected BSDE with parameters $(\tilde\psi,\tilde f,\tilde h)$, then
\begin{align}\nonumber
&\|\tilde Y^{\tilde t,\tilde x}_{\cdot\vee \tilde t}-Y^{t,x}_{\cdot\vee t}\|_{\mcS^2}^2+\|\ett_{[\tilde t,T]}\tilde Z^{\tilde t,\tilde x}-\ett_{[t,T]}Z^{t,x}\|_{\mcH^2}^2+\|\tilde K^{\tilde t,\tilde x}_{\cdot\vee \tilde t}-K^{t,x}_{\cdot\vee t}\|_{\mcS^2}^2\leq C(\|\tilde h(\cdot\vee\tilde t,X^{\tilde t,\tilde x}_{\cdot\vee \tilde t})-h(\cdot\vee t,X^{t,x}_{\cdot\vee \tilde t})\|_{\mcS^{2}}\Psi_T^{1/2}
\\
&\quad+\E\Big[|\tilde \psi(X^{\tilde t,\tilde x}_T)-\psi(X^{t,x}_T)|^{2}+\int_0^T |\ett_{[\tilde t,T]}(s)\tilde f(s,X^{\tilde t,\tilde x}_s,Y_s,Z_s)-\ett_{[t,T]}(s)f(s,X^{t,x}_s,Y_s,Z_s)|^2ds\Big]),\label{ekv:ElKaroui-diff}
\end{align}
where
\begin{align*}
\Psi_T&:=\E\Big[|\tilde \psi(X^{\tilde t,\tilde x}_T)|^{2}+|\psi(X^{t,x}_T)|^{2}+\int_{\tilde t}^T|\tilde f(s,X^{\tilde t,\tilde x}_s,0,0)|^2ds+\int_{t}^T|f(s,X^{t,x}_s,0,0)|^2ds
\\
&\quad+\sup_{s\in[\tilde t,T]}|(\tilde h(s,X^{t,x}_s))^+|^2 + \sup_{s\in[t,T]}|(h(s,X^{t,x}_s))^+|^{2}\Big].
\end{align*}
\end{thm}

\section{The case of a local driver and the corresponding impulse control problem\label{sec:impulse-control}}
In this section we consider a simplified setting in which $f$ is a function $f:[0,T]\times \R^n\times\R\times\R^d\to \R$ satisfying the requirements in the statement of Theorem~\ref{thm:ElKaroui}. We show that in this case we can relate the solution to the BSDE to an impulse control problem. We recall the system of reflected BSDEs from the introduction
\begin{align}\label{ekv:syst-bsde-simp}
\begin{cases}
  Y^{t,x}_s=\psi(X^{t,x}_T)+\int_s^T f(r,X^{t,x}_r,Y^{t,x}_r,Z^{t,x}_r)dr-\int_s^T Z^{t,x}_r dW_r+ K^{t,x}_T-K^{t,x}_s,\quad\forall s\in[t,T], \\
  Y^{t,x}_s\geq \mcM \bar Y(s,X^{t,x}_s),\:\forall s\in[t,T] \quad{\rm and}\quad\int_t^T(Y^{t,x}_s-\mcM \bar Y(s,X^{t,x}_s))dK^{t,x}_s=0.
\end{cases}
\end{align}
For $(t,x)\in [0,T]\times \R^n$ and $u\in\mcU_t$ we let the \cadlag process $X^{t,x;u}$ solve the impulsively controlled SDE
\begin{align}
X^{t,x;u}_s&=x+\int_t^s a(r,X^{t,x;u}_r)dr+\int_t^s\sigma(r,X^{t,x;u}_r)dW_r\label{ekv:forward-sde1}
\end{align}
for $s\in [t,\tau_{1})$ and
\begin{align}
X^{t,x;u}_{s}&=\Gamma(\tau_j,X^{t,x;[u]_{j-1}}_{\tau_j},\beta_j)+\int_{\tau_j}^s a(r,X^{t,x;u}_r)dr+\int_{\tau_j}^s\sigma(r,X^{t,x;u}_r)dW_r,\label{ekv:forward-sde2}
\end{align}
whenever $s\in [\tau_{j},\tau_{j+1})$ for $j=1,\ldots,N-1$ and $s\in [\tau_N,T]$ when $j=N$. Moreover, we let the pair $(P^{t,x;u},Q^{t,x;u})\in\mcS^2_{\text{cag},t}\times\mcH^2_t$ be the unique solution to the non-standard BSDE
\begin{align}
P^{t,x;u}_s&=\psi(X^{t,x;u}_T)+\int_s^T f(r,X^{t,x;u}_r,P^{t,x;u}_r,Q^{t,x;u}_r)dr-\int_s^T Q^{t,x;u}_r dW_r - \Xi^{t,x;u}_{T}+\Xi^{t,x;u}_s, \label{ekv:non-ref-bsde-simp}
\end{align}
where the impulse cost process $\Xi$ is defined as
\begin{align}\label{ekv:Xi-def}
\Xi^{t,x;u}_s:=\sum_{j=1}^N\ett_{[\tau_j<s]}\ell(\tau_j,X^{t,x;[u]_{j-1}}_{\tau_{j}},\beta_j).
\end{align}

We then introduce the following impulse control problem:\\

\textbf{Problem 1.} For $(t,x)\in[0,T]\times \R^n$, find $u^*\in\mcU_t$ such that
\begin{align*}
  P^{t,x;u^*}_t=\sup_{u\in\mcU_t} P^{t,x;u}_t.
\end{align*}
\bigskip

The main result of this section is the following:
\begin{thm}\label{thm:simp}
There exists a unique solution $(Y,Z,K)\in\bigS^2\times\bigH^2\times\bigS^2_i$ to \eqref{ekv:syst-bsde-simp}. Moreover, there is a $v\in\Pi^g_c$ such that $v(t,x)=Y^{t,x}_t$, $\Prob$-a.s., and $v$ is the unique viscosity solution in $\Pi^g$ to
\begin{align}\label{ekv:var-ineq-simp}
\begin{cases}
  \min\{v(t,x)-\mcM v(t,x),-v_t(t,x)-\mcL v(t,x)-f(t,x,v(t,x),\sigma^\top(t,x)\nabla_x v(t,x))\}=0,\\
  \quad\forall (t,x)\in[0,T)\times \R^d \\
  v(T,x)=\psi(x).
\end{cases}
\end{align}
Finally, we have the representation
\begin{align}\label{ekv:rep-simp}
  Y^{t,x}_t=\esssup_{u\in\mcU_t}P_t^{t,x;u}=P_t^{t,x;u^*}.
\end{align}
where $u^*=(\tau^*_j,\beta^*_j)_{j=1}^{N^*}\in\mcU_t$ is defined as:
\begin{itemize}
  \item $\tau^*_{j}:=\inf \big\{s \geq \tau^*_{j-1}:\:v(s,X^{t,x;[u^*]_{j-1}}_{s})=\mcM v(s,X^{t,x;[u^*]_{j-1}}_s)\big\}\wedge T$,
  \item $\beta^*_j\in\mathop{\arg\max}_{b\in U}\{v(\tau^*_{j},\Gamma(\tau^*_{j},X^{t,x;[u^*]_{j-1}}_{\tau^*_{j}},b))-\ell(\tau^*_{j},X^{t,x;[u^*]_{j-1}}_{\tau^*_{j}},b)\}$
\end{itemize}
and $N^*=\sup\{j:\tau^*_j<T\}$, with $\tau_0^*:=t$.
\end{thm}

Since the uniqueness of solutions to \eqref{ekv:var-ineq-simp} is rather standard we postpone the proof to the appendix. The proof of the existence part is based on an approximation scheme given in the following subsection.

\subsection{An approximation scheme}
The existence part of the proof of Theorem~\ref{thm:simp} is based on an approximation routine where we restrict the number of allowed interventions. We defined the following sequence of reflected BSDEs
\begin{align}\label{ekv:syst-bsde-k}
\begin{cases}
  Y^{t,x,k}_s=\psi(X^{t,x}_T)+\int_s^T f(r,X^{t,x}_r,Y^{t,x,k}_r,Z^{t,x,k}_r)dr-\int_s^T Z^{t,x,k}_r dW_r+ K^{t,x,k}_T-K^{t,x,k}_s,\quad\forall s\in[t,T], \\
  Y^{t,x,k}_s\geq \mcM \bar Y^{k-1}(s,X^{t,x}_s),\:\forall s\in[t,T] \quad{\rm and}\quad\int_t^T(Y^{t,x,k}_s-\mcM \bar Y^{k-1}(s,X^{t,x}_s))dK^{t,x,k}_s=0
\end{cases}
\end{align}
and
\begin{align}\label{ekv:syst-bsde-0}
  Y^{t,x,0}_s=\psi(X^{t,x}_T)+\int_s^T f(r,X^{t,x}_r,Y^{t,x,0}_r,Z^{t,x,0}_r)dr-\int_s^T Z^{t,x,0}_r dW_r,\quad\forall s\in[t,T].
\end{align}
We will relate the solution of these BSDEs to a sequence of impulse control problems, defined for $k\geq 1$:\\

\textbf{Problem 1.$k$} For $(t,x)\in[0,T]\times \R^n$, find $u^*\in\mcU^k_t$ such that
\begin{align*}
  P^{t,x;u^*}_t=\sup_{u\in\mcU^k_s} P^{t,x;u}_t.
\end{align*}
\bigskip

We also make use of the following lemma which is given without proof as it follows immediately from the definitions:
\begin{lem}\label{lem:mcM-monotone}
Let $u,v:[0,T]\times \R^n\to\R$ be locally bounded functions. $\mcM$ is monotone (if $u\leq v$ pointwise, then $\mcM u\leq \mcM v$). Moreover, $\mcM(u_*)$ (resp. $\mcM(u^*)$) is l.s.c.~(resp. u.s.c.).
\end{lem}
In particular, it follows that $\mcM v$ is jointly continuous whenever $v$ is.\\

Using Theorem~\ref{thm:ElKaroui} and the previous lemma we get the following intermediate result:

\begin{prop}\label{prop:simp-k}
The system \eqref{ekv:syst-bsde-k} admits a unique solution and there is a $v_k\in\Pi^g_c$ such that $v_k(t,x)=Y^{t,x}_t$ and
\begin{align}\label{ekv:bsde-opt-rep}
  v_k(t,x)=\sup_{u\in\mcU^k_t}P^{t,x;u}_t
\end{align}
for each $k\geq 0$. Moreover, an optimal control $u^{*}\in\mcU^k_t$ for \eqref{ekv:bsde-opt-rep} exists.
\end{prop}

\noindent \emph{Proof.} First note that classically \eqref{ekv:syst-bsde-0} admits a unique solution $(Y^{t,x,0},Z^{t,x,0})$ and that there is a $v_0\in\Pi^g_c$ such that $Y^{t,x,0}_s=v_0(s,X^{t,x}_s)$, $\Prob$-a.s.~for all $s\in[t,T]$. Then, since $(t,x)\mapsto \mcM g(t,x)\in\Pi^g_c$ whenever $g\in\Pi^g_c$ by Lemma~\ref{lem:mcM-monotone} it follows that $\mcM v_0\in\Pi^g_c$. On the other hand, by definition
\begin{align*}
  \mcM \bar Y^{0}(s,X^{t,x}_s)=\mcM v_0(s,X^{t,x}_s),\quad\forall s\in[t,T],\quad\Prob\text{-a.s.}
\end{align*}
and we conclude by Theorem~\ref{thm:ElKaroui} that \eqref{ekv:syst-bsde-k} admits a unique solution for $k=1$ with $v_1:=\bar Y^1\in\Pi^g_c$. Repeating this argument gives that \eqref{ekv:syst-bsde-k} admits a unique solution for arbitrary $k$ with $v_k:=\bar Y^k\in\Pi^g_c$.


For the representation, we pick $(t,x)\in[0,T)\times\R^n$ and define $u^*=(\tau^*_j,\beta^*_j)_{j=1}^{N^*}\in\mcU^k_t$ as:
\begin{itemize}
  \item $\tau^*_{j}:=\inf \big\{s \geq \tau^*_{j-1}:\:v_{k+1-j}(s,X^{t,x;[u^*]_{j-1}}_{s})=\mcM v_{k-j}(s,X^{t,x;[u^*]_{j-1}}_s)\big\}\wedge T$,
  \item $\beta^*_j\in\mathop{\arg\max}_{b\in U}\{v_{k-j}(\tau^*_{j},\Gamma(\tau^*_{j},X^{t,x;[u^*]_{j-1}}_{\tau^*_{j}},b))-\ell(\tau^*_{j},X^{t,x;[u^*]_{j-1}}_{\tau^*_{j}},b)\}$
\end{itemize}
for $j=1,\ldots,k$, and $N^*:=\max\{j\in \{0,\ldots,k\}:\tau^*_j<T\}$, with $\tau_0^*:=t$. Then, $K^{t,x}_{\tau^*_1}-K^{t,x}_{t}=0$, $\Prob$-a.s., and we find that
\begin{align*}
  Y^{t,x,k}_s&= \ett_{[\tau^*_1<T]}\mcM v_{k-1}(\tau^*_1,X^{t,x}_{\tau_1^*}) +\ett_{[\tau^*_1=T]}\psi(X^{t,x}_T) + \int_s^{\tau_1^*} f(r,X^{t,x}_r,Y^{t,x,k}_r,Z^{t,x,k}_r)dr - \int_s^{\tau_1^*}Z^{t,x,k}_r dW_r
  \\
  &= \ett_{[\tau^*_1<T]}\{v_{k-1}(\tau^*_1,X^{t,x;[u^*]_1}_{\tau_1^*})-\ell({\tau_1^*},X^{t,x}_{\tau_1^*},\beta_1^*)\} +\ett_{[\tau^*_1=T]}\psi(X^{t,x}_T)
  \\
  &\quad +\int_s^{\tau_1^*} f(r,X^{t,x}_r,Y^{t,x,k}_r,Z^{t,x,k}_r)dr - \int_s^{\tau_1^*}Z^{t,x,k}_r dW_r
\end{align*}
for all $s\in [t,\tau^*_1]$. To simplify notation later we let $(\mcY^0,\mcZ^0):=(Y^{t,x},Z^{t,x})$. Similarly, there is a triple $(\mcY^1,\mcZ^1,\mcK^1)\in\mcS^2_t\times\mcH^2_t\times\mcS^2_{t,i}$ that solves the reflected bsde
\begin{align*}
\begin{cases}
  \mcY^{1}_s=\psi(X^{t,x;[u^*]_1}_T)+\int_s^T f\big(r,X^{t,x;[u^*]_1}_r,\mcY^1,\mcZ^1_r\big)dr-\int_s^T \mcZ^1_r dW_r+ \mcK^1_T-\mcK^1_s,\quad\forall s\in[\tau^*_1,T], \\
  \mcY^{1}_s\geq \mcM v_{k-2}(s,X^{t,x;[u^*]_1}_s),\:\forall s\in[\tau^*_1,T]\quad{\rm and}\quad
  \int_{\tau_1^*}^T(\mcY^{1}_s-\mcM v_{k-2}(s,X^{t,x;[u^*]_1}_s))d\mcK^{1}_s=0
\end{cases}
\end{align*}
and continuity together with an approximation of $X^{t,x;[u^*]_1}$ and the stability result for reflected BSDEs in Proposition 3.6 of \cite{ElKaroui1} implies that $\mcY^{1}_{s}=v_{k-1}(s,X^{t,x;[u^*]_1}_{s})$, $\Prob$-a.s., for each $s\in [\tau^*_1,T]$. In particular, the proof of Proposition 2.3 in~\cite{ElKaroui1} now gives that $\mcK^1_{\tau^*_2}-\mcK^1_{\tau^*_1}=0$, $\Prob$-a.s.~and we conclude that
\begin{align*}
  \mcY^{1}_s&= \ett_{[\tau^*_2<T]}\{v_{k-2}(\tau^*_2,X^{t,x;[u^*]_2}_{\tau_2^*})-\ell({\tau_2^*},X^{t,x;[u^*]_1}_{\tau_2^*},\beta_2^*)\} +\ett_{[\tau^*_2=T]}\psi(X^{t,x;u^*}_T)
  \\
  &\quad \int_s^{\tau_2^*} f(r,X^{t,x}_r,\mcY^{1}_r,\mcZ^{1}_r)dr - \int_s^{\tau_2^*}\mcZ^{1}_r dW_r
\end{align*}
for all $s\in [\tau^*_1,\tau^*_2]$.

Repeating this process $k$ times we find that there is a sequence $(\mcY^j,\mcZ^j)_{j=0}^k\subset\mcS^2\times\mcH^2$ such that $\mcY:=\ett_{[t,\tau^*_{1}]}\mcY^0+\sum_{j=1}^k\ett_{(\tau^*_{j},\tau^*_{j+1}]}\mcY^j$ and $\mcZ:=\sum_{j=0}^k\ett_{(\tau^*_{j},\tau^*_{j+1}]}\mcZ^j$ satisfies
\begin{align*}
  \mcY_s&=\psi(X^{t,x;u^*}_T)+\int_s^{T} f\big(r,X^{t,x;u^*}_r,\mcY_r,\mcZ_r\big)dr-\int_s^{T}\mcZ_rdW_r -\sum_{j=1}^{N^*}\ett_{[s\leq \tau^*_j]}\ell({\tau_j^*},X^{t,x;[u^*]_{j-1}}_{\tau_j^*},\beta_j^*),\quad\forall s\in[t,T].
\end{align*}
By uniqueness of solutions to \eqref{ekv:non-ref-bsde-simp} we thus conclude that $Y^{t,x,k}_t=P^{t,x;u^*}_t$. Now, suppose that $\hat u\in \mcU^k_t$ is another impulse control, then
\begin{align*}
  Y^{t,x,k}_s&= \ett_{[\hat\tau_1<T]}\{v_{k-1}(\hat\tau_1,\Gamma(\hat\tau_1,X^{t,x}_{\hat\tau_1})) -\ell({\hat\tau_1},X^{t,x}_{\hat\tau_1},\hat\beta_1)\} +\ett_{[\hat\tau_1=T]}\psi(X^{t,x}_T)
  \\
  &\quad +\int_s^{\hat\tau_1} f(r,X^{t,x}_r,Y^{t,x,k}_r,Z^{t,x,k}_r)dr - \int_s^{\hat\tau_1}Z^{t,x,k}_r dW_r+K^{t,x,k}_{\hat\tau_1}-K^{t,x,k}_s
\end{align*}
for all $s\in [t,\hat\tau_1]$. Arguing as above gives that there is a sequence $(\hat\mcY^j,\hat\mcZ^j,\hat\mcK^j)_{j=0}^k\subset\mcS^2\times\mcH^2\times\mcS^2_i$ such that letting $\hat\mcY:=\ett_{[t,\hat\tau_{1}]}\hat\mcY^0+\sum_{j=0}^k\ett_{(\hat\tau_{j},\hat\tau_{j+1}]}\hat\mcY^j$, $\mcZ:=\sum_{j=0}^k\ett_{(\hat\tau_{j},\hat\tau_{j+1}]}\hat\mcZ^j$ and $\hat\mcK_s:=\sum_{j=0}^{k}\ett_{[\hat \tau_j<s]}\{\hat\mcK^j_{s\wedge\hat\tau_{j+1}}-\hat\mcK^j_{\hat\tau_j}\}$, with $\hat\tau_0:=t$, implies that $(\hat\mcY,\hat\mcZ,\hat\mcK)\subset\mcS^2_{\text{cag},t}\times\mcH^2_t\times\mcS^2_{i,t}$ satisfies
\begin{align*}
  \hat\mcY_s&= \psi(X^{t,x;\hat u}_T)+\int_s^{T} f\big(r,X^{t,x;\hat u}_r,\hat\mcY_r,\hat\mcZ_r\big)dr-\int_s^{T} \hat\mcZ_rdW_r-\sum_{j=1}^{\hat N}\ett_{[s\leq\hat\tau_j]}\ell({\hat \tau_j},X^{t,x;[\hat u]_{j-1}}_{\hat \tau_j},\hat \beta_j)
  \\
  &\quad+\hat\mcK_{T}-\hat \mcK_{s},\quad\forall s\in[t,T]
\end{align*}
and $\hat\mcY_t=Y^{t,x}_t$. Now, comparison gives that $P^{t,x;\hat u}_t\leq\hat\mcY_t$ and we conclude that $P^{t,x;\hat u}_t \leq Y^{t,x,k}_t$.\qed\\

Before we proceed to give the proof of Theorem~\ref{thm:simp} we need some preliminary estimates, which we give in the following subsection.

\subsection{Some preliminary estimates}
\begin{prop}\label{prop:SDEmoment}
For each $p\geq 1$, there is a $C>0$ such that
\begin{align}\label{ekv:SDEmoment}
\E\Big[\sup_{s\in[\zeta,T]}|X^{t,x;u}_s|^{p}\Big|\mcF^t_\zeta\Big]\leq C(1+|X^{t,x;u}_\zeta|^{p}),
\end{align}
$\Prob$-a.s.~for all $(t,x)\in [0,T]\times\R^n$ and $(\zeta,u)\in [t,T]\times\mcU_t$.
\end{prop}

\noindent\emph{Proof.} We use the shorthand $X^j:=X^{t,x;[u]_j}$. By Assumption~\ref{ass:onSFDE}.(\ref{ass:onSFDE-Gamma}) we get for $s\in [\tau_{j},T]$, using integration by parts, that
\begin{align*}
|X^{j}_s|^2&= |X^{j}_{\zeta\vee\tau_{j}}|^2+2\int_{(\zeta\vee\tau_{j})+}^s X^{j}_{r}dX^{j}_r+\int_{(\zeta\vee\tau_{j})+}^s d[X^{j},X^{j}]_r
\\
&\leq K^2_\Gamma\vee |X^{{j-1}}_{\zeta\vee\tau_{j}}|^2+2\int_{(\zeta\vee\tau_{j})+}^s X^{j}_{r} dX^{j}_r+\int_{(\zeta\vee\tau_{j})+}^s d[X^{j},X^{j}]_r.
\end{align*}
Now, either $|X^{{j-1}}_{\tau_{j}}|\leq K_\Gamma$ in which case
\begin{align*}
|X^{j}_s|^2\leq |X^{j}_{\zeta}|^2\vee K^2_\Gamma+2\int_{(\zeta\vee\tau_{j})+}^s X^{j}_{r} dX^{j}_r+\int_{(\zeta\vee\tau_{j})+}^s d[X^{j},X^{j}]_r.
\end{align*}
or $|X^{{j-1}}_{\tau_{j}}|> K_\Gamma$ implying that
\begin{align*}
|X^{j}_s|^2&\leq K^2_\Gamma\vee |X^{{j-2}}_{\zeta\vee\tau_{j-1}}|^2+2\int_{(\zeta\vee\tau_{j-1})+}^{\tau_{j}} X^{j-1}_{r} dX^{j-1}_r+\int_{(\zeta\vee\tau_{j-1})+}^{\tau_{j}} d[X^{j-1},X^{j-1}]_r
\\
&\quad+2\int_{(\zeta\vee\tau_{j})+}^s X^{j}_{r} dX^{j}_r+\int_{(\zeta\vee\tau_{j})+}^s d[X^{j},X^{j}]_r.
\end{align*}
In the latter case the same argument can be repeated and we conclude that
\begin{align}\label{ekv:X2-bound}
|X^{j}_s|^2&\leq |X^{j}_{\zeta}|^2\vee K_\Gamma^2+\sum_{i=j_0}^{j} \Big\{2\int_{(\zeta\vee\tilde\tau_{i})+}^{s\wedge\tilde\tau_{i+1}} X^{i}_{r}dX^{i}_r+\int_{(\zeta\vee\tilde\tau_{i})+}^{s\wedge\tilde\tau_{i+1}} d[X^{i},X^{i}]_r\Big\},
\end{align}
where $\tilde\tau_0=-1$, $\tilde\tau_i=\tau_i$ for $i=1,\ldots,j$ and $\tilde\tau_{j+1}=\infty$ and $j_0:=\max\{i\in \{1,\ldots,j\}:|X^{{i-1}}_{\tau_{i}}|\leq K_\Gamma\}\vee 0$.

Now, since $X^{i}$ and $X^{j}$ coincide on $[0,\tau_{i+1\wedge j+1})$ we have
\begin{align*}
\sum_{i=j_0}^{j}\int_{(\zeta\vee\tilde\tau_{i})+}^{s\wedge\tilde\tau_{i+1}} X^{i}_{r} dX^{i}_r
&=\int_{\zeta\vee\tau_{j_0}}^s X^{j}_{r}a(r,X^{j}_r)dr+\int_{\zeta\vee\tau_{j_0}}^{s}X^{j}_{r}\sigma(r,X^{j}_r)dW_r,
\end{align*}
and
\begin{align*}
\sum_{i=j_0}^{j} \int_{(\zeta\vee\tilde\tau_{i})+}^{s\wedge\tilde\tau_{i+1}} d[X^{i},X^{i}]_r&=\int_{\zeta\vee\tau_{j_0}}^{s} \sigma^2(r,X^{j}_r)dr.
\end{align*}
Inserted in \eqref{ekv:X2-bound} this that
\begin{align}\nonumber
|X^{j}_s|^2&\leq |X^{j}_{\zeta}|^2\vee K_\Gamma^2+\int_{\tau_{j_0}}^s (2X^{j}_{s}a(r,X^{j}_r)+\sigma^2(r,X^{j}_r))dr+2\int_{\tau_{j_0}}^{s}X^{j}_{r}\sigma(r,X^{j}_r)dW_r
\\
&\leq |X^{j}_{\zeta}|^2+C\Big(1+\int_{\zeta}^{s}|X^{j}_{r}|^2dr + \sup_{\eta\in[\zeta,s]}\Big|\int_{\zeta}^{\eta}X^{j}_r\sigma(r,X^{j}_r)dW_r\Big|\Big)\label{ekv:X-bound1}
\end{align}
for all $s\in [\zeta,T]$. The Burkholder-Davis-Gundy inequality now gives that for $p\geq 2$,
\begin{align*}
\E\Big[\sup_{r\in[\zeta,s]}|X^{j}_r|^{p}\Big|\mcF^t_\zeta \Big]\leq |X^{j}_{\zeta}|^2+C\big(1+\E\Big[\int_{\zeta}^{s}|X^{i}_{r}|^{p}dr+\big(\int_{\zeta}^{s}|X^{j}_r|^4 dr\big)^{p/4}\Big|\mcF^t_\zeta\Big]\big)
\end{align*}
and Gr\"onwall's lemma gives that for $p\geq 4$,
\begin{align*}
\E\Big[\sup_{s\in[\zeta,T]}|X^{j}_s|^{p}\Big|\mcF^t_\zeta\Big]&\leq C(1+ |X^{j}_{\zeta}|^{p}),
\end{align*}
$\Prob$-a.s., where the constant $C=C(T,p)$ does not depend on $u$ or $j$ and \eqref{ekv:SDEmoment} follows by letting $j\to\infty$ on both sides and using Fatou's lemma. The result for general $p\geq 1$ follows by Jensen's inequality.\qed\\

For $(t,x)\in[0,T]\times\R^n$ and $u\in\mcU_t$ we let $(\check P^{t,x;u},\check Q^{t,x;u})$ be the unique solution to the following standard BSDE
\begin{align}\label{ekv:bsde-trad}
  \check P_s^{t,x;u}=\psi(X^{t,x;u}_T)+\int_s^Tf(r,X^{t,x;u}_r,\check P^{t,x;u}_r,\check Q^{t,x;u}_r)dr-\int_s^T \check Q^{t,x;u}_rdW_r.
\end{align}
Combining classical results (see \eg \cite{ElKaroui2}) with Proposition~\ref{prop:SDEmoment}, we have
\begin{align}\nonumber
  &\E\Big[\sup_{s\in [t,T]}|\check P_s^{t,x;u}|^2+\int_t^T|\check Q^{t,x;u}_s|^2ds\Big]
  \\
  &\leq C\E\Big[|\psi(X^{t,x;u}_T)|^2+\int_t^T|f(r,X^{t,x;u}_r,0,0)|^2dr\Big]\leq C(1+|x|^{2\rho}),\label{ekv:bsde-trad-moment}
\end{align}
for all $u\in\mcU_t$.

Using the comparison principle we easily deduce the following moment estimates:
\begin{prop}\label{prop:BSDEmoment}
We have,
\begin{align}\label{ekv:Ybnd}
  |\sup_{u\in\mcU_t}P^{t,x;u}_t|\leq C(1+|x|^\rho)
\end{align}
and for each $k\geq 0$, there is a $C>0$ such that
\begin{align}\label{ekv:BSDEmoment}
\E\Big[\sup_{s\in[t,T]}|P^{t,x;u}_s|^{2}+\int_t^T| Q^{t,x;u}_s|^2ds\Big]\leq C(1+|x|^{2\rho}),
\end{align}
$\Prob$-a.s.~for all $(t,x)\in [0,T]\times\R^n$ and $u\in\mcU^k_t$.
\end{prop}

\noindent\emph{Proof.} The first statement follows by repeated application of the comparison principle which gives that $\check P^{t,x;\emptyset}_t\leq\sup_{u\in\mcU_t}P^{t,x;u}_t \leq \sup_{u\in\mcU_t}\check P^{t,x;u}_t$ and using~\eqref{ekv:bsde-trad-moment}.

The second statement follows by noting that for fixed $k\geq 0$, there is a $C>0$ such that
\begin{align*}
\E[|\Xi^{t,x;u}_{T}|^2]&\leq C(1+\E[\sup_{s\in[t,T]}|X^{t,x;u}_{s}|^{2\rho}])\leq C(1+|x|^{2\rho})
\end{align*}
for all $u\in \mcU^k_t$.\qed\\

In the following lemma we use the above estimates to derive a bound on the expected number of interventions in an optimal control for Problem 1.$k$. This bound plays an important role in the convergence analysis employed later on.

\begin{lem}\label{lem:EN-bound}
There is a $C>0$ such that whenever $u^*$ is an optimal control to Problem 1.$k$ for some $k\geq 0$, then $\E[N^*]\leq C(1+|x|^\rho)$.
\end{lem}

\noindent \emph{Proof.} To simplify notation we let $(X,P,Q):=(X^{t,x;u^*},P^{t,x;u^*},Q^{t,x;u^*})$ and $X^j=X^{t,x;[u^*]_j}$ and get that
\begin{align*}
P_s&=\psi(X_T)+\int_s^{T}f(r,X_r,P_r,Q_r)dr-\int_s^{T} Q_rdW_r-\sum_{\tau_j\geq s}\ell(\tau_j,X^{j-1}_{\tau_j},\beta_j).
\end{align*}
Letting
\begin{align*}
 \zeta_1(s):=\frac{f(s,X_s,P_s,Q_s) -f(s,X_s,0,Q_s)}{P_s}\ett_{[P_s\neq 0]}
\end{align*}
and
\begin{align*}
 \zeta_2(s):=\frac{f(s,X_s,0,Q_s) - f(s,X_s,0,0)}{|Q_s|^2}(Q_s)^\top
\end{align*}
we have by the Lipschitz continuity of $f$ that $|\zeta_1(s)|\vee|\zeta_2(s)|\leq k_f$. Using Ito's formula we find that
\begin{align*}
P_s&=R_{s,T}\psi(X_T)+\int_s^{T}R_{s,r}f(r,X_r,0,0)dr -\int_s^TR_{s,r}Q_rdW_r-\sum_{j=1}^{N^*} R_{s,\tau^*_{j}}\ett_{[\tau^*_j\geq s]}\ell(\tau^*_j,X^{j-1}_{\tau^*_{j}},\beta^*_j)
\end{align*}
with $R_{s,r}:=e^{\int_s^{r}(\zeta_1(v)-\frac{1}{2}|\zeta_2(v)|^2)dv+\frac{1}{2}\int_s^{r}\zeta_2(v)dW_v}$. Since the intervention costs are positive, taking the conditional expectation on both sides and using Proposition~\ref{prop:SDEmoment} gives
\begin{align*}
P_s&\leq\E\Big[R_{s,T}\psi(X_T)+\int_s^{T}R_{s,r}f(r,X_r,0,0)dr \Big|\mcF^t_s\Big]
\\
&\leq C\Big(1+\E\big[R_{s,T}^2\big|\mcF^t_s\big]^{1/2}\E\Big[\sup_{r\in [s,T]}|X_r|^{2\rho} \Big|\mcF^t_s\Big]^{1/2}\Big)
\\
&\leq C(1+|X_s|^\rho)
\end{align*}
On the other hand, since $u^*$ is an optimal control,
\begin{align*}
  P_s\geq v_0(s,X_s)\geq -C(1+|X_s|^\rho),
\end{align*}
$\Prob$-a.s., for some $C>0$ (independent of $(t,s,x)$ and $k$). Proposition~\ref{prop:SDEmoment} then gives
\begin{align*}
\E\Big[\sup_{s\in[t,T]}|P_s|^2\Big]&\leq C(1+|x|^{2\rho}),
\end{align*}
where $C>0$ does not depend on $k$. Next, we derive a bound on the $\mcH^2_t$-norm of $Q$. Applying Ito's formula to $|P_s|^{2}$ we get
\begin{align}\nonumber
|P_t|^{2}+\int_t^T| Q_s|^2ds&=\psi^2(X_T)+2\int_t^T P_sf(s,X_s,P_s,Q_s)ds-2\int_t^T P_sQ_sdW_s
\\
&\quad-\sum_{j=1}^{N^*}(2P^{j-1}_{\tau^*_{j}}\ell(\tau^*_j,X^{j-1}_{\tau^*_{j}},\beta^*_j) + \ell^2(\tau^*_j,X^{j-1}_{\tau^*_{j}},\beta^*_j)),\label{ekv:from-ito}
\end{align}
where $P^{j-1}$ is $P$ without the $j-1$ first intervention costs. Since the intervention costs are nonnegative, we have
\begin{align*}
-\sum_{j=1}^{N^*} (2P^{j-1}_{\tau^*_{j}}\ell(\tau^*_j,X^{j-1}_{\tau^*_{j}},\beta^*_j)+\ell^2(\tau^*_j,X^{j-1}_{\tau^*_{j}},\beta^*_j)) &\leq 2\sup_{s\in [t,T]}|P_{s}|\sum_{j=1}^{N^*} \ell(\tau^*_j,X^{j-1}_{\tau^*_{j}},\beta^*_j)
\\
&\leq \kappa \sup_{s\in [t,T]}|P_{s}|^2+\frac{1}{\kappa}\Big(\sum_{j=1}^{N^*} \ell(\tau^*_j,X^{j-1}_{\tau^*_{j}},\beta^*_j)\Big)^2
\end{align*}
for any $\kappa>0$. Inserted in \eqref{ekv:from-ito} and using the Lipschitz property of $f$ this gives
\begin{align}\nonumber
|P_t|^{2}+\int_t^T| Q_s|^2ds&\leq \psi^2(X_T)+(C+\kappa)\sup_{s\in[t,T]}|P_s|^2+\int_t^T(|f(s,X_s,0,0)|^2+\frac{1}{2}|Q_s|^2)ds
\\
&\quad -2\int_t^T P_sQ_sdW_s+\frac{1}{\kappa}\Big(\sum_{j=1}^{N^*} \ell(\tau^*_j,X^{j-1}_{\tau^*_{j}},\beta^*_j)\Big)^2.\label{ekv:PQ-bound}
\end{align}
Now, as $u^*\in\mcU^k$, it follows that the stochastic integral is uniformly integrable and thus a martingale. To see this, note that the Burkholder-Davis-Gundy inequality gives
\begin{align*}
\E\Big[\sup_{s\in[t,T]}\Big|\int_t^s P_rQ_rdW_r\Big|\Big]\leq C\E\Big[\Big(\int_t^T |P_sQ_s|^2ds\Big)^{1/2}\Big]\leq C\E\Big[\sup_{s\in [t,T]}|P_s|^2+\int_t^T |Q_s|^2ds\Big]
\end{align*}
where the right-hand side is finite by \eqref{ekv:BSDEmoment}. Taking expectations on both sides of \eqref{ekv:PQ-bound} thus gives
\begin{align*}
\E\Big[\int_t^T| Q_s|^2ds\Big]&\leq C(1+\kappa)(1+|x|^{2\rho})+\frac{2}{\kappa}\E\Big[\Big(\sum_{j=1}^{N^*} \ell(\tau^*_j,X^{j-1}_{\tau^*_{j}},\beta^*_j)\Big)^2\Big].
\end{align*}
Finally,
\begin{align*}
\E[N^*]&\leq \frac{1}{\delta}\E\Big[\Big(\sum_{j=1}^{N^*} \ell(\tau^*_j,X^{j-1}_{\tau^*_{j}},\beta^*_j)\Big)^2\Big]^{1/2}
\end{align*}
and
\begin{align*}
\E\Big[\Big(\sum_{j=1}^{N^*} \ell(\tau^*_j,X^{j-1}_{\tau^*_{j}},\beta^*_j)\Big)^2\Big]&\leq C\E\Big[|P_t|^2+|\psi(X_T)|^2+\int_t^{T}|f(r,X_r,P_r,Q_r)|^2dr+\int_t^{T}|Q_r|^2dr\Big]
\\
&\leq C\E\Big[|\psi(X_T)|^2+\sup_{s\in[t,T]}|P_s|^2+\int_t^{T}(|f(s,X_s,0,0)|^2+|Q_s|^2)ds\Big]
\\
&\leq C(1+\kappa)(1+|x|^{2\rho})+\frac{C}{\kappa}\E\Big[\Big(\sum_{j=1}^{N^*} \ell(\tau^*_j,X^{j-1}_{\tau^*_{j}},\beta^*_j)\Big)^2\Big]
\end{align*}
and the lemma follows by choosing $\kappa$ sufficiently large.\qed\\

\begin{lem}\label{lem:v_k-unif-conv}
There is a $C>0$ such that
\begin{align}\label{ekv:v_k-unif-conv}
  |v_k(t,x)-v_{k-1}(t,x)|\leq \frac{C}{k}(1+|x|^{2\rho})
\end{align}
\end{lem}

\noindent\emph{Proof.} We let $u^*\in\mcU^k_t$ be an optimal control for Problem 1.$k$. To simplify notation we let $(P,Q)=(P^{t,x;u^*},Q^{t,x;u^*})$ and set\\ $(\hat P,\hat Q)=(P^{t,x;[u^*]_{k-1}},Q^{t,x;[u^*]_{k-1}})$. Then, with $X:=X^{t,x;u^*}$ and $\hat X:=X^{t,x;[u^*]_{k-1}}$, we have $\hat X_s=X_s$ for all $s\in[0,\tau^*_{k})\cap[0,T]$. This gives
\begin{align*}
P_t-\hat P_t&=\psi(X_T)-\psi(\hat X_T) +\int_t^{T}(f(s,X_s,P_s,Q_s)- f(s,\hat X_s,\hat P_s,\hat Q_s))ds
\\
&\quad -\int_t^{T}(Q_s- \hat Q_s)dW_s+\Xi^{t,x;[u^{*}]_{k-1}}_{T}-\Xi^{t,x;u^{*}}_{T}
\\
&\leq \ett_{[N^{*}=k]}\Big(R_{t,T}(\psi(X_T)-\psi(\hat X_T)) +\int_t^{T}R_{t,s}(f(s,X_s,P_s,Q_s)- f(s,\hat X_s,P_s,Q_s))ds\Big)
\\
&\quad-\int_t^{T}R_{t,s}(Q_s- \hat Q_s)dW_s
\end{align*}
for some $R_{t,s}:=e^{\int_t^{s}(\zeta_1-\frac{1}{2}|\zeta_2(r)|^2)dr+\frac{1}{2}\int_t^{s}\zeta_2(r)dW_r}$, with $|\zeta_1(r)|\vee|\zeta_2(r)|\leq k_f$. Taking expectation on both sides and using the Cauchy-Schwartz inequality gives
\begin{align*}
\E[P_t-\hat P_t]&\leq\E\Big[\ett_{[N^{*}=k]}\Big(R_{t,T}(\psi(X_T) - \psi(\hat X_T))+\int_t^{T}R_{t,s}(f(s,X_s,P_s,Q_s)- f(s,\hat X_s,P_s,Q_s))ds\Big) \Big]
\\
&\leq C(1+|x|^{\rho})\E\big[\ett_{[N^{*}=k]}\big]^{1/2}.
\end{align*}
Now, Lemma~\ref{lem:EN-bound} implies that
\begin{align*}
\E\big[\ett_{[N^{*}=k]}\big]\leq  \frac{C(1+|x|^\rho)}{k}
\end{align*}
from which \eqref{ekv:v_k-unif-conv} follows.\qed\\

\subsection{Proof of Theorem~\ref{thm:simp}}
\noindent\emph{Proof of Theorem~\ref{thm:simp}} Since the sequence $(v_k)_{k\geq 0}$ is non-decreasing and uniformly bounded by a polynomial, there is a $v\in\Pi^g$ such that $v_k\nearrow v$ as $k\to\infty$. Moreover, Lemma~\ref{ekv:v_k-unif-conv} implies that the limit is jointly continuous, thus $v\in\Pi^g_c$.

Now, by Theorem~\ref{thm:ElKaroui} there is, for each $(t,x)\in [0,T]\times\R^n$, a unique triple $(Y^{t,x},Z^{t,x},K^{t,x})\in\mcS_t^2\times\mcH^2_t\times\mcS^2_{i,t}$ such that
\begin{align}\label{ekv:rbsde-simp-2}
\begin{cases}
  Y^{t,x}_s=\psi(X^{t,x}_T)+\int_s^T f(r,X^{t,x}_r,Y^{t,x}_r,Z^{t,x}_r)dr-\int_s^T Z^{t,x}_r dW_r+ K^{t,x}_T-K^{t,x}_s,\quad\forall s\in[t,T], \\
  Y^{t,x}_s\geq \mcM v(X^{t,x}_s),\quad\forall s\in[t,T] \quad{\rm and}\quad\int_t^T(Y^{t,x}_s-\mcM v(X^{t,x}_s))dK^{t,x}_s=0.
\end{cases}
\end{align}
On the other hand, we easily deduce by \eqref{ekv:ElKaroui-diff} that $Y^{t,x}=\lim_{k\to\infty}Y^{t,x,k}$ pointwisely, and we conclude that $Y^{t,x}_t=v(t,x)$, $\Prob$-a.s. In particular, this implies that $(Y,Z,K)\in \bigS^2\times\bigH^2\times\bigS^2_i$ solves \eqref{ekv:syst-bsde-simp}.

Concerning the representation, we get by repeating the proof of Proposition~\ref{prop:simp-k} that, whenever $(Y,Z,K)\in\bigS^2\times\bigH^2\times\bigS^2_i$ satisfies \eqref{ekv:syst-bsde-simp}, there is a unique pair $(\mcY,\mcZ)\in\mcS^2_{\text{cag},t}\times\mcH^2_t$ such that $\mcY_t=Y^{t,x}_t$, $\Prob$-a.s., and
\begin{align*}
  \mcY_s&=\int_s^{\tau^*_k} f\big(r,X^{t,x;u^*}_r,\mcY_r,\mcZ\big)dr -\int_s^{\tau^*_k} \mcZ_rdW_r -\sum_{j=1}^{k\wedge N^*}\ett_{[s\leq\tau^*_j]}\ell({\tau_j^*},X^{t,x;[u^*]_{j-1}}_{\tau_j^*},\beta_j^*)
  \\
  &\quad+\ett_{[\tau^*_k=T]}\psi(X^{t,x;u^*}_T)+\ett_{[\tau^*_k<T]}v(\tau^*_{k},X^{t,x;[u^*]_{k}}_{\tau_k^*}),
\end{align*}
where $u^*$ is now the impulse control in the statement of the theorem. The sequence $v(\tau^*_{k},X^{t,x;[u^*]_{k}}_{\tau_k^*})$ is uniformly bounded in $L^2(\Prob)$ and we conclude that $u^*\in\mcU_t$. In particular, this implies that $N^*$ is $\Prob$-a.s.~finite. Taking the limit as $k\to\infty$ and using that \eqref{ekv:non-ref-bsde-simp} admits a unique solution, we conclude that $Y^{t,x}_t=P^{t,x;u^*}_t$. Repeating the comparison part in the proof of Proposition~\ref{prop:simp-k} then gives that $u^*$ is an optimal control for Problem 1. Since this holds for any $(Y,Z,K)\in \bigS^2_c\times\bigH^2\times\bigS^2_i$ that solves \eqref{ekv:syst-bsde-simp}, uniqueness of solutions to \eqref{ekv:syst-bsde-simp} follows.

Moreover, from Theorem~\ref{thm:ElKaroui} it immediately follows that $v$ solves~\eqref{ekv:var-ineq-simp} and by the comparison result for viscosity solutions to \eqref{ekv:var-ineq-simp} in Proposition~\ref{app:prop:comp-visc} the solution is unique.\qed\\

\begin{rem}
Letting $\mcU^t_s$ be the subset of $\mcU_t$ with $\tau_1\geq s$, $\Prob$-a.s., we may extend the representation to $Y^{t,x}_s=\esssup_{u\in\mcU^t_s}P^{t,x;u}_s$ for all $s\in[t,T]$.
\end{rem}

\section{The general setting\label{sec:general}}
We now turn to the general setting of a non-local driver. Existence will again follow by an approximation routine and we introduce the following sequence of systems of BSDEs
\begin{align}\label{ekv:rbsde_k}
\begin{cases}
  Y^{t,x,k}_s=\psi(X^{t,x}_T)+\int_s^T f(r,X^{t,x}_r,\bar Y^{k-1}(r,\cdot),Z^{t,x,k}_r)dr-\int_s^T Z^{t,x,k}_rdW_r+ K^{t,x,k}_T-K^{t,x,k}_s,\:\forall s\in[t,T] \\
  Y^{t,x,k}_s\geq \mcM \bar Y^{k}(s,X^{t,x}_s),\:\forall s\in[t,T]\quad\text{and}\quad \int_0^T(Y^{t,x,k}_s-\mcM \bar Y^{k}(s,X^{t,x}_s))dK^{t,x,k}_s=0.
\end{cases}
\end{align}
for $k\geq 1$, with $\bar Y^{0}\equiv 0$.

\begin{prop}
There is a sequence $((Y^{k},Z^{k},K^{k})\in \bigS^2\times\bigH^2\times\bigS^2_i)_{k\geq 0}$ that satisfies the recursion in \eqref{ekv:rbsde_k}.
\end{prop}

\emph{Proof.} We need to show that for each $k\geq 1$, there is a $v_{k-1}\in\Pi^g_c$ such that $Y^{t,x,k-1}_t=v_{k-1}(t,x)$ for all $(t,x)\in [0,T]\times\R^n$. However, for $k=1$ this is immediate by the definition. Now, the result follows by using Theorem~\ref{thm:simp} and induction.\qed\\

For $(t,\gamma)\in[0,T]\times\R_+$ and $\alpha\in\mcA_t$, we let $(\Psi^{t,\gamma;\alpha},\Theta^{t,\gamma;\alpha})\in\mcS^2_t\times\mcS^2_{i,t}$ (with $\Theta^{t,\gamma;\alpha}_t=0$) solve the one-dimensional reflected SDE
\begin{align}\label{ekv:rsde}
\begin{cases}
  \Psi^{t,\gamma;\alpha}_s=\gamma^2\vee K_\Gamma^2+(4C_{a,\sigma}+2C_{a,\sigma}^2)\int_t^s(1+\Psi^{t,\gamma;\alpha}_r)dr+4C_{a,\sigma}\int_t^s (1+\Psi^{t,\gamma;\alpha}_r)\alpha_r dW_r+ \Theta^{t,\gamma;\alpha}_s \\
  \Psi^{t,\gamma;\alpha}_s\geq \gamma^2\vee K_\Gamma^2\text{ and } \int_t^T(\Psi^{t,\gamma;\alpha}_s-(\gamma^2\vee K_\Gamma^2))d\Theta^{t,\gamma;\alpha}_s=0,
\end{cases}
\end{align}
For $(t,\gamma)\in[0,T]\times \R_+$, we then set $R^{t,\gamma;\alpha}_s:=\sqrt{\Psi^{t,\gamma;\alpha}_s}$ and note that classically, we have
\begin{align*}
\E\Big[\sup_{s\in [t,T]}|R^{t,\gamma;\alpha}_s|^p\Big]\leq C(1+|\gamma\vee K_\Gamma|^p),
\end{align*}
for all $p\geq 2$.

\begin{lem}\label{lem:R-bounds-X}
For each $(t,x)\in [0,T]\times \R^d$ and $u\in\mcU_t$, there is an $\alpha\in\mcA_t$ such that $|X^{t,x,u}_s|\leq R^{t,|x|;\alpha}_s$ for all $s\in [t,T]$, $\Prob$-a.s.
\end{lem}

\noindent\emph{Proof.} Since
\begin{align*}
|2x a(r,x)+\sigma^2(r,x)|\leq (4C_{a,\sigma}+2C_{a,\sigma}^2)(1+|x|^2),
\end{align*}
it follows from \eqref{ekv:X-bound1} that we can always choose $\alpha\in\mcA_t$ such that
\begin{align*}
2C_{a,\sigma}(1+\Psi^{t,\gamma,\xi;\alpha}_r)\alpha_r=X^{j}_{r}\sigma(r,X^{j}_r)
\end{align*}
and the statement holds by \eqref{ekv:X-bound1}.\qed\\

For $\varphi\in\Pi^g_c$ we let $(Y^{\varphi},Z^\varphi,K^\varphi)\in \bigS^2_c\times\bigH^2\times\bigS^2_i$ be the unique solution to
\begin{align}\label{ekv:seq-bsde-loc-2}
\begin{cases}
  Y^{t,x,\varphi}_s=\psi(X^{t,x}_T)+\int_s^T f(r,X^{t,x}_r,\varphi(r,\cdot),Z^{t,x,\varphi}_r)dr-\int_s^T Z^{t,x,\varphi}_r dW_r+ K^{t,x,\varphi}_T-K^{t,x,\varphi}_s,\quad\forall s\in[t,T],\\
  Y^{t,x,\varphi}_s\geq \mcM \bar Y^{\varphi}(s,X^{t,x}_s),\:\forall s\in[t,T]\quad\text{and}\quad \int_0^T(Y^{t,x,\varphi}_s-\mcM \bar Y^{\varphi}(s,X^{t,x}_s))dK^{t,x,\varphi}_s=0
\end{cases}
\end{align}
and note that letting $(P^{t,x,\varphi;u},Q^{t,x,\varphi;u})\in\mcS^2_{\text{cag},t}\times\mcH^2_t$ solve
\begin{align}
P^{t,x,\varphi;u}_s&=\psi(X^{t,x;u}_T)+\int_s^T f(r,X^{t,x;u}_r,\varphi(r,\cdot),Q^{t,x,\varphi;u}_r)dr-\int_s^T Q^{t,x,\varphi;u}_r dW_r - \Xi^{t,x;u}_{T}+\Xi^{t,x;u}_s, \label{ekv:non-ref-bsde}
\end{align}
Theorem~\ref{thm:simp} gives that
\begin{align*}
  Y^{t,x,\varphi}_t=\sup_{u\in\mcU_t}P^{t,x,\varphi;u}_t.
\end{align*}

\begin{prop}\label{prop:contraction}
There is a $\kappa>0$ such that for all $\varphi,\tilde \varphi\in \Pi^g_c$ and $\gamma>0$, we have
\begin{align}\label{ekv:int-contraction}
  \sup_{\alpha\in\mcA}\E\Big[\int_0^Te^{\kappa t}\sup_{x\in\Lambda_f(R^{0,\gamma;\alpha}_t)}|\bar Y^{\tilde \varphi}(t,x)-\bar Y^{\varphi}(t,x)|^2dt\Big]\leq \frac{1}{2}\sup_{\alpha\in\mcA}\E\Big[\int_0^T e^{\kappa t}\sup_{x\in\Lambda_f(R^{0,\gamma;\alpha}_t)}|\tilde \varphi(t,x)-\varphi(t,x)|^2dt\Big].
\end{align}
Furthermore, there is a $C>0$ such that
\begin{align}\label{ekv:sup-contraction}
  \sup_{t\in[0,T]}\sup_{x\in\Lambda_f(\gamma)}|\bar Y^{\tilde \varphi}(t,x)-\bar Y^{\varphi}(t,x)|^2\leq C\sup_{\alpha\in\mcA}\E\Big[\int_0^T \sup_{x\in\Lambda_f(R^{0,\gamma;\alpha}_t)}|\tilde \varphi(t,x)-\varphi(t,x)|^2dt\Big]
\end{align}
for each $\gamma>0$.
\end{prop}

\noindent\emph{Proof.} Let $u^*,\tilde u^*\in\mcU_t$ be optimal strategies for $P^{t,x,\varphi;u}$ and $P^{t,x,\tilde \varphi;u}$, respectively, so that
\begin{align*}
  Y^{t,x,\varphi}_t=P^{t,x,\varphi;u^*}_t \quad \text{and} \quad Y^{t,x,\tilde \varphi}_t=P^{t,x,\tilde \varphi;\tilde u^*}_t.
\end{align*}
Then
\begin{align*}
  Y^{t,x,\varphi}_t- Y^{t,x,\tilde \varphi}_t&=P^{t,x,\tilde \varphi;u^*}_t-P^{t,x,\tilde \varphi;\tilde u^*}_t
  \\
  &\leq P^{t,x,\tilde \varphi;u^*}_t-P^{t,x,\tilde \varphi;u^*}_t.
\end{align*}
Since a similar inequality holds in the opposite direction, we find that
\begin{align}\label{ekv:YbyP}
  |Y^{t,x,\varphi}_t- Y^{t,x,\tilde \varphi}_t|&\leq \sup_{u\in \mcU_t}|P^{t,x,\tilde \varphi;u}_t-P^{t,x,\tilde \varphi;u}_t|.
\end{align}
For $u\in\mcU_t$, let $(P,Q):=(P^{t,x,\varphi;u},Q^{t,x,\varphi;u})$ and $(\tilde P,\tilde Q):=(P^{t,x,\tilde \varphi;u},Q^{t,x,\tilde \varphi;u})$ and note that for $\kappa>0$, It\^o's formula applied to $e^{\kappa \cdot}|\tilde P-P|^2$ gives
\begin{align*}
  e^{\kappa t}|\tilde P_t-P_t|^2&+\int_t^Te^{\kappa s}|\tilde Q_s-Q_s|^2ds=-2\int_t^Te^{\kappa s}(\tilde P_s-P_s)(\tilde Q_s-Q_s)dW_s-\kappa \int_t^Te^{\kappa s}|\tilde P_s-P_s|^2ds
  \\
  &\quad + 2\int_t^Te^{\kappa s}(\tilde P_s-P_s)(f(s,X^{t,x;u}_s,\tilde \varphi(s,\cdot),\tilde Q_s)-f(s,X^{t,x;u}_s,\varphi(s,\cdot),Q_s))ds.
\end{align*}
By assumption
\begin{align*}
  |f(s,X^{t,x;u}_s,\varphi(s,\cdot),Q_s)-f(s,X^{t,x;u}_s,\tilde \varphi(s,\cdot),\tilde Q_s)|\leq k_f(\sup_{x'\in \Lambda_f(|X^{t,x;u}_s|)}|\tilde \varphi(s,x')-\varphi(s,x')| + |\tilde Q_s-Q_s|).
\end{align*}
Hence, taking the expectation and using inequalities $2Cxy\leq (Cx)^2+y^2$ and $2xy\leq x^2/\sqrt{\kappa}+\sqrt{\kappa}y^2$ gives
\begin{align*}
  e^{\kappa t}|\tilde P_t-P_t|^2&\leq (C^2+C\sqrt{\kappa}-\kappa)\E\Big[ \int_t^Te^{\kappa s}|\tilde P_s-P_s|^2ds \Big]
  \\
  &\quad+\frac{C}{\sqrt{\kappa}}\E\Big[\int_t^Te^{\kappa s}\sup_{x'\in \Lambda_f(|X^{t,x;u}_s|)}|\tilde \varphi(s,x')-\varphi(s,x')|^2ds \Big].
\end{align*}
Now, pick $\kappa_0>0$ such that $\kappa_0\geq C^2+C\sqrt{\kappa_0}$ and note that for each $\kappa\geq\kappa_0$, we have
\begin{align*}
  e^{\kappa t}|P_t-\tilde P_t|^2&\leq \frac{C}{\sqrt{\kappa}}\E\Big[\int_t^Te^{\kappa s}\sup_{x'\in \Lambda_f(|X^{t,x;u}_s|)}|\tilde \varphi(s,x')-\varphi(s,x')|^2ds \Big]
  \\
  &\leq \frac{C}{\sqrt{\kappa}}\sup_{\alpha\in\mcA_t}\E\Big[\int_t^Te^{\kappa s}\sup_{x'\in \Lambda_f(|R^{t,|x|;\alpha}_s|)}|\tilde \varphi(s,x')-\varphi(s,x')|^2ds \Big],
\end{align*}
where the last inequality follows from Lemma~\ref{lem:R-bounds-X}. Since the right-hand side is non-decreasing in $|x|$ and independent of $u$, \eqref{ekv:YbyP} now gives that
\begin{align}\label{ekv:contraction-intermediate}
  e^{\kappa t}\sup_{x\in\Lambda_f(\gamma)}|\bar Y^{\tilde \varphi}(t,x)-\bar Y^{\varphi}(t,x)|^2&\leq \frac{C}{\sqrt{\kappa}}\sup_{\alpha\in\mcA_t}\E\Big[\int_t^Te^{\kappa s}\sup_{x'\in \Lambda_f(R^{t,\gamma;\alpha}_s)}|\tilde \varphi(s,x')-\varphi(s,x')|^2ds \Big],
\end{align}
for any $\gamma\geq 0$. In particular, as both sides are continuous in $\gamma$ a standard dynamic programming argument gives that for any $\alpha_1\in \mcA$, we have
\begin{align*}
  \E\Big[e^{\kappa t}\sup_{x\in\Lambda_f(R^{0,\gamma;\alpha_1}_t)}|\bar Y^{\tilde \varphi}(t,x)-\bar Y^{\varphi}(t,x)|^2\Big]&\leq \frac{C}{\sqrt{\kappa}}\sup_{\alpha\in\mcA}\E\Big[\int_t^Te^{\kappa s}\sup_{x'\in \Lambda_f(R^{0,\gamma;\alpha_1\oplus_t\alpha}_s)}|\tilde \varphi(s,x')-\varphi(s,x')|^2ds \Big].
\end{align*}
Taking the supremum with respect to $\alpha_1$ on the right hand side and once again relying on a standard dynamic programming argument gives that
\begin{align*}
  \E\Big[e^{\kappa t}\sup_{x\in\Lambda_f(R^{0,\gamma;\alpha_1}_t)}|\bar Y^{\tilde \varphi}(t,x)-\bar Y^{\varphi}(t,x)|^2\Big]&\leq \frac{C}{\sqrt{\kappa}}\sup_{\alpha\in\mcA}\E\Big[\int_0^Te^{\kappa s}\sup_{x'\in \Lambda_f(R^{0,\gamma;\alpha}_s)}|\tilde \varphi(s,x')-\varphi(s,x')|^2ds\Big].
\end{align*}
Integrating with respect to time and using Fubini's theorem, we find that
\begin{align*}
  \E\Big[\int_0^T e^{\kappa t}\sup_{x\in\Lambda_f(R^{0,\gamma;\alpha_1}_t)}|\bar Y^{\tilde \varphi}(t,x)-\bar Y^{\varphi}(t,x)|^2 dt\big]&\leq \frac{CT}{\sqrt{\kappa}}\sup_{\alpha\in\mcA}\E\Big[\int_0^Te^{\kappa s}\sup_{x'\in \Lambda_f(R^{0,\gamma;\alpha}_s)}|\tilde \varphi(s,x')-\varphi(s,x')|^2ds\Big]
\end{align*}
after which taking the supremum with respect to $\alpha_1\in\mcA$ and choosing $\kappa\geq 2(CT)^2\vee\kappa_0$ gives the first inequality. To get \eqref{ekv:sup-contraction} we note that comparison gives that $R^{0,\gamma;\alpha}_s\geq R^{t,\gamma;\alpha}_s$ for all $s\in [t,T]$ and $\alpha\in\mcA$. From \eqref{ekv:contraction-intermediate} we thus get that
\begin{align*}
  \sup_{x\in\Lambda_f(\gamma)}|\bar Y^{\tilde \varphi}(t,x)-\bar Y^{\varphi}(t,x)|^2&\leq C\sup_{\alpha\in\mcA}\E\Big[\int_0^T\sup_{x'\in \Lambda_f(R^{0,\gamma;\alpha}_s)}|\tilde \varphi(s,x')-\varphi(s,x')|^2ds \Big]
\end{align*}
from which \eqref{ekv:sup-contraction} is immediate since the right hand side is independent of $t$.\qed\\

We now introduce the norm $\|\cdot\|_{\gamma}$ on the space of jointly continuous functions of polynomial growth, $\Pi^g_c$, defined as
\begin{align*}
  \|\varphi\|^2_{\gamma}:=\sup_{\alpha\in\mcA}\E\Big[\int_0^Te^{\kappa t}\sup_{x\in \Lambda_f(R^{0,\gamma;\alpha}_t)}|\varphi(t,x)|^2dt\Big],
\end{align*}
with $\kappa>0$ as in Proposition~\ref{prop:contraction} and note that under $\|\cdot\|_{\gamma}$, the map $\Phi:\Pi^g_c\to\Pi^g_c$ that maps $\varphi$ to $\bar Y^{\varphi}$ is a contraction.

\begin{cor}\label{cor:unif-bound}
There are constants $C>0$ and $p\geq 0$ such that $|\bar Y^{k}(t,x)|\leq C(1+|x|^p)$ for all $(t,x)\in[0,T]\times\R^d$ and all $k\geq 0$.
\end{cor}

\noindent\emph{Proof.} First, we note that \eqref{ekv:int-contraction} and the triangle inequality implies that
\begin{align*}
  \|\bar Y^{k}\|_{\gamma}&\leq  \|\bar Y^{k}-\bar Y^{k-1}\|_{\gamma}+\|\bar Y^{k-1}\|_{\gamma}\leq \frac{1}{2}\|\bar Y^{k-1}-\bar Y^{k-2}\|_{\gamma}+\|\bar Y^{k-1}\|_{\gamma}\leq \frac{1}{2^{k-1}}\|\bar Y^{1}-\bar Y^{0}\|_{\gamma}+\|\bar Y^{k-1}\|_{\gamma}.
\end{align*}
However, as a similar scheme holds for $\|\bar Y^{k-1}\|_{\gamma}$ and since $\bar Y^{0}\equiv 0$ we conclude that
\begin{align*}
  \|\bar Y^{k}\|_{\gamma}&\leq  \sum_{j=1}^k\frac{1}{2^{j-1}}\|\bar Y^{1}\|_{\gamma}\leq 2\|\bar Y^{1}\|_{\gamma}.
\end{align*}
On the other hand, as $\bar Y^{1}\in\Pi^g_c$ there are constants $C>0$ and $p\geq 2$ such that $|\bar Y^{1}(t,x)|\leq C(1+|x|^p)$ and we conclude that
\begin{align*}
  \|\bar Y^{1}\|^2_{\gamma}&=\sup_{\alpha\in\mcA}\E\Big[\int_0^Te^{\kappa t}\sup_{x\in \Lambda_f(R^{0,\gamma;\alpha}_t)}|\bar Y^{1}(t,x)|^2dt\Big]
  \\
  &\leq C\Big(1+\sup_{\alpha\in\mcA}\E\Big[\sup_{t\in [0,T]}|R^{0,\gamma;\alpha}_t|^{2p}\Big]\Big)
  \\
  &\leq C(1+|\gamma|^{2p})
\end{align*}
implying the existence of a $C>0$ such that $\|\bar Y^{k}\|_{\gamma}\leq C(1+|\gamma|^{p})$ for all $k\geq 0$. Now, \eqref{ekv:sup-contraction} gives that
\begin{align*}
  \sup_{t\in[0,T]}\sup_{x\in\Lambda_f(\gamma)}|\bar Y^{k}(t,x)-\bar Y^{1}(t,x)|^2&\leq C\sup_{\alpha\in\mcA}\E\Big[\int_0^T \sup_{x\in\Lambda_f(R^{0,\gamma;\alpha}_t)}|\bar Y^{k-1}(t,x)|^2dt\Big]
  \\
  &\leq C(1+|\gamma|^{2p})
\end{align*}
where the constants $C>0$ and $p\geq 2$ do not depend on $k$ and the desired bound follows.\qed\\

Letting $v_k(t,x):=\bar Y^{k}(t,x)$, Proposition~\ref{prop:contraction} and Corollary~\ref{cor:unif-bound} implies that there is a $v\in \Pi^g$ such that for each $\gamma>0$ we have $\|v_k-v\|_{\gamma}\to 0$ as $k\to\infty$.

\begin{thm}
$v$ is the unique viscosity solution in $\Pi^g_c$ to \eqref{ekv:var-ineq}.
\end{thm}

\noindent \emph{Proof.} First, \eqref{ekv:sup-contraction} implies that the convergence is uniform on compact subsets of $[0,T]\times\R^n$ and since $v_k$ is jointly continuous for each $k\geq 0$ we conclude that $v$ is also jointly continuous. This in turn gives that $\Phi(v)$ is well defined and we conclude that $\Phi(v)=v$ establishing existence of a solution to \eqref{ekv:syst-bsde-new}. Moreover, if $(\tilde Y,\tilde Z,\tilde K)$ is another solution, then $\tilde v(t,x):=\tilde Y^{t,x}_t$ must also satisfy $\Phi(\tilde v)=\tilde v$. However, then repeated use of the contraction property in \eqref{ekv:int-contraction} gives that $\|\tilde v-v\|_{\gamma,\kappa}=0$ and by continuity we conclude that $\tilde v=v$ implying by uniqueness of solutions to \eqref{ekv:syst-bsde-simp} as obtained in Theorem~\ref{thm:simp} that \eqref{ekv:syst-bsde-new} admits a unique solution.

Utilizing, once more, the connection between reflected BSDEs and obstacle problems we conclude that $v$ solves \eqref{ekv:var-ineq}. Suppose now that there exists another function $\tilde v\in\Pi^g_c$ that solves \eqref{ekv:var-ineq} and let $\bar v=\Phi(\tilde v)$, then by Theorem~\ref{thm:simp} we conclude that $\bar v\in\Pi^g_c$ is the unique solution to
\begin{align*}
\begin{cases}
  \min\{\bar v(t,x)-\mcM \bar v(t,x),-\bar v_t(t,x)-\mcL \bar v(t,x)-f(t,x,\tilde v(t,\cdot),\sigma^\top(t,x)\nabla_x \bar v(t,x))\}=0,\\
  \quad\forall (t,x)\in[0,T)\times \R^d \\
  \bar v(T,x)=\psi(x),
\end{cases}
\end{align*}
But then $\bar v=\tilde v$ implying that $\tilde v$ is a fixed point of $\Phi$ and since $v$ is the only fixed point of $\Phi$ in the set of jointly continuous functions of polynomial growth we conclude that $\tilde v=v$.\qed\\

\begin{cor}
$(Y^v,Z^v,K^v)$ is the unique solution in $\bigS^2\times\bigH^2\times\bigS^2_i$ to \eqref{ekv:syst-bsde-new}.
\end{cor}

\appendix
\section{Uniqueness of viscosity solutions in the local framework\label{app:uni-simp}}
By Theorem~\ref{thm:simp}, there is $v\in\Pi^g_c$ that solves (in viscosity sense) the quasi-variational inequality
\begin{align}\label{app:ekv:var-ineq-simp}
\begin{cases}
  \min\{v(t,x)-\mcM v(t,x),-v_t(t,x)-\mcL v(t,x)-f(t,x,v(t,x),\sigma^\top(t,x)\nabla_x v(t,x))\}=0,\\
  \quad\forall (t,x)\in[0,T)\times \R^d \\
  v(T,x)=\psi(x).
\end{cases}
\end{align}
In this section we show that $v$ is the only viscosity solution in $\Pi^g$ to \eqref{app:ekv:var-ineq-simp}. We need the following lemma:

\begin{lem}\label{lem:is-super}
Let $v$ be a supersolution to \eqref{ekv:var-ineq} satisfying
\begin{align*}
\forall (t,x)\in[0,T]\times\R^d,\quad |v(t,x)|\leq C(1+|x|^{2\varrho})
\end{align*}
for some $\varrho>0$. Then there is a $\gamma_0 > 0$ such that for any $\gamma>\gamma_0$ and $\theta > 0$, the function $v + \theta e^{-\gamma t}(1+((|x|-K_\Gamma)^+)^{2\varrho+2})$ is also a supersolution to \eqref{ekv:var-ineq}.
\end{lem}

\noindent \emph{Proof.} With $w:=v + \theta e^{-\gamma t}(1+((|x|-K_\Gamma)^+)^{2\varrho+2})$ we note that, since $v$ is a supersolution and $\theta e^{-\gamma T}(1+((|x|-K_\Gamma)^+)^{2\varrho+2})\geq 0$, we have $w(T,x)\geq v(T,x)\geq \psi(x)$ so that the terminal condition holds. Moreover, we have
\begin{align*}
&w(t,x)-\sup_{b\in U}\{w(t,\Gamma(t,x,b))-\ell(t,x,b)\}
\\
&=v(t,x) + \theta e^{-\gamma t}(1+((|x|-K_\Gamma)^+)^{2\varrho+2})
\\
&\quad-\sup_{b\in U}\{v(t,\Gamma(t,x,b)) + \theta e^{-\gamma t}(1+((|\Gamma(t,x,b)|-K_\Gamma)^+)^{2\varrho+2})-\ell(t,x,b)\}
\\
&\geq v(t,x)- \sup_{b\in U}\{v(t,\Gamma(t,x,b))-\ell(t,x,b)\}
\\
&\quad+\theta e^{-\gamma t}\{1+((|x|-K_\Gamma)^+)^{2\varrho+2}- \sup_{b\in U}(1+((|\Gamma(t,x,b)|-K_\Gamma)^+)^{2\varrho+2})\}.
\end{align*}
Since $v$ is a supersolution, we have
\begin{align*}
  v(t,x)- \sup_{b\in U}\{v(t,\Gamma(t,x,b))-\ell(t,x,b)\}\geq 0
\end{align*}
Now, either $|x|\leq K_\Gamma$ in which case it follows by \eqref{ekv:imp-bound} that $|\Gamma(t,x,b)|\leq K_\Gamma$ or $|x|> K_\Gamma$ and \eqref{ekv:imp-bound} gives that $|\Gamma(t,x,b)|\leq |x|$. We conclude that
\begin{align*}
  w(t,x)- \sup_{b\in U}\{w(t,\Gamma(t,x,b))-\ell(t,x,b)\}\geq 0.
\end{align*}
Next, let $\varphi\in C^{1,2}([0,T]\times\R^d\to\R)$ be such that $\varphi-w$ has a local maximum of 0 at $(t_0,x_0)$ with $t_0<T$. Then $\tilde \varphi(t,x):=\varphi (t,x)-\theta e^{-\gamma t}(1+((|x|-K_\Gamma)^+)^{2\varrho+2})\in C^{1,2}([0,T]\times\R^d\to\R)$ and $\tilde \varphi-v$ has a local maximum of 0 at $(t_0,x_0)$. Since $v$ is a viscosity supersolution, we have
\begin{align*}
    &-\partial_t(\varphi(t,x)-\theta e^{-\gamma t}(1+((|x|-K_\Gamma)^+)^{2\varrho+2}))-\mcL (\varphi(t,x)-\theta e^{-\gamma t}(1+((|x|-K_\Gamma)^+)^{2\varrho+2}))
    \\
    &-f(t,x,\varphi(t,x)-\theta e^{-\gamma t}(1+((|x|-K_\Gamma)^+)^{2\varrho+2}),\sigma^\top(t,x)\nabla_x (\varphi(t,x)-\theta e^{-\gamma t}(1+((|x|-K_\Gamma)^+)^{2\varrho+2})))\geq 0.
\end{align*}
Consequently,
\begin{align*}
&-\partial_t\varphi(t,x)-\mcL \varphi(t,x)-f(t,x,\varphi(t,x),\sigma^\top(t,x)\nabla_x \varphi(t,x))
\\
&\geq \theta \gamma e^{-\gamma t}(1+((|x|-K_\Gamma)^+)^{2\varrho+2})-\theta\mcL e^{-\gamma t}(1+((|x|-K_\Gamma)^+)^{2\varrho+2})
\\
&\quad f(t,x,\varphi(t,x)-\theta e^{-\gamma t}(1+((|x|-K_\Gamma)^+)^{2\varrho+2}),\sigma^\top(t,x)\nabla_x (\varphi(t,x)-\theta e^{-\gamma t}(1+((|x|-K_\Gamma)^+)^{2\varrho+2})))
\\
&\quad-f(t,x,\varphi(t,x),\sigma^\top(t,x)\nabla_x \varphi(t,x))
\\
&\geq \theta \gamma e^{-\gamma t}(1+((|x|-K_\Gamma)^+)^{2\varrho+2})-\theta C(1+\varrho) e^{-\gamma t}(1+((|x|-K_\Gamma)^+)^{2\varrho+2})
\\
&\quad -k_f(1+\varrho)\theta e^{-\gamma t}(1+((|x|-K_\Gamma)^+)^{2\varrho+2}),
\end{align*}
where the right hand side is non-negative for all $\theta> 0$ and all $\gamma>\gamma_0$ for some $\gamma_0>0$.\qed\\

We have the following result, the proof of which we omit since it is classical:
\begin{lem}\label{lem:integ-factor}
For any $\lambda\in\R$, a locally bounded function $v:[0,T]\times \R^d\to\R$ is a viscosity supersolution (resp. subsolution) to \eqref{ekv:var-ineq} if and only if $\tilde v(t,x):=e^{\lambda t}v(t,x)$ is a viscosity supersolution (resp. subsolution) to
\begin{align}\label{ekv:var-ineq-if}
\begin{cases}
  \min\{\tilde v(t,x)-\sup_{b\in U}\{\tilde v(t,\Gamma(t,x,b))-e^{\lambda t}\ell(t,x,b)\},-\tilde v_t(t,x)+\lambda \tilde v(t,x)-\mcL \tilde v(t,x)\\
  -e^{\lambda t}f(t,x,e^{-\lambda t}\tilde v(t,x),e^{-\lambda t}\sigma^\top(t,x)\nabla_x \tilde v(t,x))\}=0,\quad\forall (t,x)\in[0,T)\times \R^d \\
  \tilde v(T,x)=e^{\lambda T}\psi(x).
\end{cases}
\end{align}
\end{lem}
\begin{rem}
Here, it is important to note that $\tilde \ell(t,x):=e^{\lambda t}\ell(t,x)$, $\tilde f(t,x,y,z):=-\lambda y+e^{\lambda t}f(t,x,e^{-\lambda t}y,e^{-\lambda t}z)$ and $\tilde \psi(x):=e^{\lambda T}\psi(x)$ satisfy Assumption~\ref{ass:on-coeff}. In particular, this implies that Lemma~\ref{lem:is-super} holds for supersolutions to \eqref{ekv:var-ineq-if} as well.
\end{rem}

We have the following comparison result for viscosity solutions in $\Pi^g$:

\begin{prop}\label{app:prop:comp-visc}
Let $v$ (resp. $u$) be a supersolution (resp. subsolution) to \eqref{ekv:var-ineq}. If $u,v\in \Pi^g$, then $u\leq v$.
\end{prop}

\noindent \emph{Proof.} First, we note that it is sufficient to show that the statement holds for solutions to~\eqref{ekv:var-ineq-if} for some $\lambda\in \R$. We thus assume that $v$ (resp.~$u$) is a viscosity supersolution (resp.~subsolution) to \eqref{ekv:var-ineq-if} for $\lambda\in \R$ specified below. Furthermore, we may without loss of generality assume that $v$ is l.s.c.~and $u$ is u.s.c.

By assumption, $u,v\in \Pi^g$, which implies that there are $C>0$ and $\varrho>0$ such that
\begin{align}\label{ekv:uv-bound}
|v(t,x)|+|u(t,x)|\leq C(1+|x|^{2\varrho}).
\end{align}
Now, for any $\gamma>0$ we only need to show that
\begin{align*}
w(t,x)&=w^{\theta,\gamma}(t,x):=v(t,x)+\theta e^{-\gamma t}(1+((|x|-K_\Gamma)^+)^{2\varrho+2})
\\
&\geq u(t,x)
\end{align*}
for all $(t,x)\in[0,T]\times\R^d$ and any $\theta>0$. Then the result follows by taking the limit $\theta\to 0$. We know from Lemma~\ref{lem:is-super} that there is a $\gamma_0>0$ such that $w$ is a supersolution to \eqref{ekv:var-ineq-if} for each $\gamma\geq\gamma_0$ and $\theta>0$. We thus assume that $\gamma\geq\gamma_0$.

We search for a contradiction and assume that there is a $(t_0,x_0)\in [0,T]\times \R^d$ such that $u(t_0,x_0)>w(t_0,x_0)$. By \eqref{ekv:uv-bound}, there is for each $\theta>0$ a $R\geq K_\Gamma$ such that
\begin{align*}
w(t,x)>u(t,x),\quad\forall (t,x)\in[0,T]\times\R^d,\:|x|>R.
\end{align*}
Our assumption thus implies that there is a point $(\bar t,\bar x)\in[0,T)\times B_R$ (the open unit ball of radius $R$ centered at 0) such that
\begin{align*}
\max_{(t,x)\in[0,T]\times\R^d}(u(t,x)-w(t,x))&=\max_{(t,x)\in[0,T)\times B_R}(u(t,x)-w(t,x))
\\
&=u(\bar t,\bar x)-w(\bar t,\bar x)=\eta>0.
\end{align*}
We first show that there is at least one point $(t^*,x^*)\in[0,T)\times B_R$ such that
\begin{enumerate}[a)]
  \item $u(t^*,x^*)-w(t^*,x^*)= \eta$ and
  \item $u(t^*,x^*)>\sup_{b\in U}\{u(t^*,\Gamma(t^*,x^*,b))-\ell(t^*,b)\}$.
\end{enumerate}
We again argue by contradiction and assume that $u(t,x)=\sup_{b\in U}\{u(t,\Gamma(t,x,b))-\ell(t,b)\}$ for all $(t,x)\in A:=\{(s,y)\in[0,T]\times\R^d: u(s,y)-w(s,y)=\eta\}$. Indeed, as $u$ is u.s.c. and $\Gamma$ is continuous, there is a $b_1$ such that
\begin{align}\label{ekv:equal}
u(\bar t,\bar x)=\sup_{b\in U}\{u(\bar t,\Gamma(\bar t,\bar x,b))-\ell(\bar t,b)\}=u(\bar t,\Gamma(\bar t,\bar x,b_1))-\ell(\bar t,b_1).
\end{align}
Now, set $x_1=\Gamma(\bar t,\bar x,b_1)$ and note that since
 \begin{align*}
|\Gamma(t,x,b)|\leq R,\quad \forall(t,x,b)\in [0,T]\times \bar B_R\times U
\end{align*}
it follows that $x_1\in \bar B_R$. Moreover, as $w$ is a supersolution it satisfies
\begin{align*}
  w(\bar t,\bar x)- (w(\bar t,\Gamma(\bar t,\bar x,b_1))-\ell(\bar t,\bar x,b_1))\geq 0
\end{align*}
or
\begin{align*}
  - w(\bar t,x_1))\geq -w(\bar t,\bar x)-\ell(t,\bar x,b_1)
\end{align*}
and we conclude from \eqref{ekv:equal} that
\begin{align*}
  u(\bar t,x_1)- w(\bar t,x_1)&\geq u(\bar t,\bar x)+\ell(\bar t,\bar x,b_1)-(w(\bar t,\bar x)+\ell(t,\bar x,b_1))
  \\
  &=u(\bar t,\bar x)-w(\bar t,\bar x)=\eta.
\end{align*}
Hence, $(\bar t,x_1)\in A$ and by our assumption it follows that there is a $b_2\in U$ such that
\begin{align*}
u(\bar t,x_1)=u(\bar t,\Gamma(\bar t,x_1,b_2))-\ell(\bar t,b_2)
\end{align*}
and a corresponding $x_2:=\Gamma(\bar t,x_1,b_2)\in B_R$. Now, this process can be repeated indefinitely to find a sequence $(x_j,b_j)_{j\geq 1}$ in $B_R\times U$ such that for any $l\geq 0$ we have
\begin{align*}
  u(\bar t,\bar x)=u(\bar t,x_l)-\sum_{j=1}^{l}\ell(\bar t,x_{j-1},b_j),
\end{align*}
with $x_0:=\bar x$. However, as $\ell\geq \delta>0$ we get a contradiction by letting $l\to\infty$ while noting that $|u(t,x)|$ is bounded on $[0,T]\times \bar B_R$. We can thus find a $(t^*,x^*)\in [0,T)\times B_R$ such that \emph{a)} and \emph{b)} above holds.

Since $f$ is Lipschitz in $y$ and $z$ for $(t,x)\in[0,T]\times\bar B_R$, the remainder of the proof follows along the lines of the proof of Proposition 4.1 in~\cite{Morlais13} and is included only for the sake of completeness.

Next, we assume without loss of generality that $\varrho\geq 2$ and define
\begin{align*}
\Phi_n(t,x,y):=u(t,x)-w(t,x)-\varphi_n(t,x,y),
\end{align*}
where
\begin{align*}
  \varphi_n(t,x,y):=\frac{n}{2}|x-y|^{2\varrho}+|x-x^*|^2+|y-\bar y|^2+(t-t^*)^2.
\end{align*}
Since $u$ is u.s.c.~and $w$ is l.s.c.~there is a triple $(t_n,x_n,y_n)\in[0,T]\times \bar B_R\times \bar B_R$ (with $\bar B_R$ the closure of $B_R$) such that
\begin{align*}
  \Phi_n(t_n,x_n,y_n)=\max_{(t,x,y)\in [0,T]\times \bar B_R\times \bar B_R}\Phi_n(t,x,y).
\end{align*}
Now, the inequality $2\Phi_n(t_n,x_n,y_n)\geq \Phi_n(t_n,x_n,x_n)+\Phi_n(t_n,y_n,y_n)$ gives
\begin{align*}
n|x_n-y_n|^{2\varrho}\leq u(t_n,x_n)-u(t_n,y_n)+w(t_n,x_n)-w(t_n,y_n).
\end{align*}
Consequently, $n|x_n-y_n|^{2\varrho}$ is bounded (since $u$ and $w$ are bounded on $[0,T]\times\bar B_R\times\bar B_R$) and $|x_n-y_n|\to 0$ as $n\to\infty$. We can, thus, extract subsequences $n_l$ such that $(t_{n_l},x_{n_l},y_{n_l})\to (\tilde t,\tilde x,\tilde x)$ as $l\to\infty$. Since
\begin{align*}
u(t^*,x^*)-w(t^*,x^*)\leq \Phi_n(t_n,x_n,y_n)\leq u(t_n,x_n)-w(t_n,y_n),
\end{align*}
it follows that
\begin{align*}
u(t^*,x^*)-w(t^*,x^*)&\leq \limsup_{l\to\infty} \{u(t_{n_l},x_{n_l})-w(t_{n_l},y_{n_l})\}
\\
&\leq u(\tilde t,\tilde x)-w(\tilde t,\tilde x)
\end{align*}
and as the righthand side is dominated by $u(t^*,x^*)-w(t^*,x^*)$ we conclude that
\begin{align*}
  u(\tilde t,\tilde x)-w(\tilde t,\tilde x)=u(t^*,x^*)-w(t^*,x^*).
\end{align*}
In particular, this gives that $\lim_{l\to\infty}\Phi_n(t_{n_l},x_{n_l},y_{n_l})=u(\tilde t,\tilde x)-w(\tilde t,\tilde x)$ which implies that
\begin{align*}
  \limsup_{l\to\infty} n_l|x_{n_l}-y_{n_l}|^{2\varrho}= 0
\end{align*}
and
\begin{align*}
  (t_{n_l},x_{n_l},y_{n_l})\to (t^*,x^*,x^*).
\end{align*}
We can thus extract a subsequence $(\tilde n_l)_{l\geq 0}$ of $(n_l)_{l\geq 0}$ such that $t_{\tilde n_l}<T$, $|x_{\tilde n_l}|<R$ and
\begin{align*}
  u(t_{\tilde n_l},x_{\tilde n_l})-w(t_{\tilde n_l},x_{\tilde n_l})\geq \frac{\eta}{2}.
\end{align*}
Moreover, since $\sup_{b\in U}\{u(t,\Gamma(t,x,b))-\tilde\ell(t,b)\}$ is u.s.c.~(see Lemma~\ref{lem:mcM-monotone}) and $u(t_{\tilde n_l},x_{\tilde n_l})\to u(t^*,x^*)$ there is an $l_0\geq 0$ such that
\begin{align*}
  u(t_{\tilde n_l},x_{\tilde n_l})-\sup_{b\in U}\{u(t_{\tilde n_l},\Gamma(t_{\tilde n_l},x_{\tilde n_l},b))-\tilde\ell(t_{\tilde n_l},b)\}>0,
\end{align*}
for all $l\geq l_0$. To simplify notation we will, from now on, denote $(\tilde n_l)_{l\geq l_0}$ simply by $n$.\\

By Theorem 8.3 of~\cite{UsersGuide} there are $(p^u_n,q^u_n,M^u_n)\in \bar J^{2,+}u(t_n,x_n)$ and $(p^w_n,q^w_n,M^w_n)\in \bar J^{2,+}w(t_n,y_n)$, where $\bar J^{2,+}$ is the limiting superjet, such that
\begin{align*}
\begin{cases}
  p^u_n-p^w_n=\partial_t\varphi_n(t_n,x_n,y_n)=2(t_n-t^*)
  \\
  q^u_n=D_x\varphi_n(t_n,x_n,y_n)=n\varrho(x-y)|x-y|^{2\varrho-2}+2(x-x^*)
  \\
  q^w_n=-D_y\varphi_n(t_n,x_n,y_n)=n\varrho(x-y)|x-y|^{2\varrho-2}+2(x-x^*)
\end{cases}
\end{align*}
and for every $\epsilon>0$,
\begin{align*}
  \left[\begin{array}{cc} M^n_x & 0 \\ 0 & -M^n_y\end{array}\right]\leq B(t_n,x_n,y_n)+\epsilon B^2(t_n,x_n,y_n),
\end{align*}
where $B(t_n,x_n,y_n):=D^2_{(x,y)}\varphi_n(t_n,x_n,y_n)$. Now, we have
\begin{align*}
  D^2_{(x,y)}\varphi_n(t,x,y)=\left[\begin{array}{cc} D_x^2\varphi_n(t,x,y) & D^2_{yx}\varphi_n(t,x,y) \\ D^2_{xy}\varphi_n(t,x,y) & D_y^2\varphi_n(t,x,y)\end{array}\right]
  = \left[\begin{array}{cc} n\xi(x,y)+2 I & -n\xi(x,y) \\ -n\xi(x,y) & n\xi(x,y) +2 I \end{array}\right]
\end{align*}
where $I$ is the identity-matrix of suitable dimension and
\begin{align*}
  \xi(x,y):=\varrho|x-y|^{2\varrho-4}\{|x-y|^2I+2(\varrho-1)(x-y)(x-y)^\top\}.
\end{align*}
In particular, since $x_n$ and $y_n$ are bounded, choosing $\epsilon:=\frac{1}{n}$ gives that
\begin{align}\label{ekv:mat-bound}
  \tilde B_n:=B(t_n,x_n,y_n)+\epsilon B^2(t_n,x_n,y_n)\leq Cn|x_n-y_n|^{2\varrho-2}\left[\begin{array}{cc} I & -I \\ -I & I \end{array}\right]+C I.
\end{align}
By the definition of viscosity supersolutions and subsolutions we have that
\begin{align*}
&-p^u_n+\lambda u(t_n,x_n)-a^\top(t_n,x_n)q^u_n-\frac{1}{2}\trace [\sigma^\top(t_n,x_n)M^u_n\sigma(t_n,x_n)]
\\
&-e^{\lambda t_n}f(t_n,x_n,e^{-\lambda t_n}u(t_n,x_n),e^{-\lambda t_n}\sigma^\top(t_n,x_n)q^u_n)\leq 0
\end{align*}
and
\begin{align*}
&-p^w_n+\lambda w(t_n,y_n)-a^\top(t_n,y_n)q^w_n-\frac{1}{2}\trace [\sigma^\top(t_n,y_n)M^w_n\sigma(t_n,y_n)]
\\
&-e^{\lambda t_n}f(t_n,y_n,e^{-\lambda t_n}w(t_n,y_n),e^{-\lambda t_n}\sigma^\top(t_n,x_n)q^w_n)\geq 0.
\end{align*}
Combined, this gives that
\begin{align*}
\lambda (u(t_n,x_n)-w(t_n,y_n))&\leq p^u_n+a^\top(t_n,x_n)q^u_n+\frac{1}{2}\trace [\sigma^\top(t_n,x_n)M^u_n\sigma(t_n,x_n)]
\\
&+e^{\lambda t_n}f(t_n,x_n,e^{-\lambda t_n}u(t_n,x_n),e^{-\lambda t_n}\sigma^\top(t_n,x_n)q^u_n)
\\
&-p^w_n-a^\top(t_n,y_n)q^w_n-\frac{1}{2}\trace [\sigma^\top(t_n,y_n)M^w_n\sigma(t_n,y_n)]
\\
&-e^{\lambda t_n}f(t_n,y_n,e^{-\lambda t_n}w(t_n,y_n),e^{-\lambda t_n}\sigma^\top(t_n,x_n)q^w_n)
\end{align*}
Collecting terms we have that
\begin{align*}
p^u_n-p^w_n&=2(t_n-t^*)
\end{align*}
and since $a$ is Lipschitz continuous in $x$ and bounded on $\bar B_R$, we have
\begin{align*}
  a^\top(t_n,x_n)q^u_n-a^\top(t_n,y_n)q^w_n&\leq  (a^\top(t_n,x_n)-a^\top(t_n,y_n))n\varrho(x_n-y_n)|x_n-y_n|^{2\varrho-2}
  \\
  &\quad+C(|x_n-x^*|+|y_n-x^*|)
  \\
  &\leq C(n|x_n-y_n|^{2\varrho}+|x_n-x^*|+|y_n-x^*|),
\end{align*}
where the right-hand side tends to 0 as $n\to\infty$. Let $s_x$ denote the $i^{\rm th}$ column of $\sigma(t_n,x_n)$ and let $s_y$ denote the $i^{\rm th}$ column of $\sigma(t_n,y_n)$ then by the Lipschitz continuity of $\sigma$ and \eqref{ekv:mat-bound}, we have
\begin{align*}
  s_x^\top M^u_n s_x-s_y^\top M^w_n s_y&=\left[\begin{array}{cc} s_x^\top  & s_y^\top \end{array}\right]\left[\begin{array}{cc} M^u_n  & 0 \\ 0 &-M^w_n\end{array}\right]\left[\begin{array}{c} s_x \\ s_y \end{array}\right]
  \\
  &\leq \left[\begin{array}{cc} s_x^\top  & s_y^\top \end{array}\right]\tilde B_n\left[\begin{array}{c} s_x \\ s_y \end{array}\right]
  \\
  &\leq C(n|x_n-y_n|^{2\varrho}+|x_n-y_n|)
\end{align*}
and we conclude that
\begin{align*}
  \limsup_{n\to\infty}\frac{1}{2}\trace [\sigma^\top(t_n,x_n)M^u_n\sigma(t_n,x_n)-\sigma^\top(t_n,y_n)M^w_n\sigma(t_n,y_n)]\leq 0.
\end{align*}
Finally, we have that
\begin{align*}
  &e^{\lambda t_n}f(t_n,x_n,e^{-\lambda t_n}u(t_n,x_n),e^{-\lambda t_n}\sigma^\top(t_n,x_n)q^u_n)-e^{\lambda t_n}f(t_n,y_n,e^{-\lambda t_n}w(t_n,y_n),e^{-\lambda t_n}\sigma^\top(t_n,x_n)q^w_n)
  \\
  &\leq k_f(u(t_n,x_n)-w(t_n,y_n)+|\sigma^\top(t_n,x_n)q^u_n-\sigma^\top(t_n,x_n)q^w_n|)
  \\
  &\quad + e^{\lambda t_n}|f(t_n,x_n,e^{-\lambda t_n}u(t_n,x_n),e^{-\lambda t_n}\sigma^\top(t_n,x_n)q^u_n)-f(t_n,y_n,e^{-\lambda t_n}u(t_n,x_n),e^{-\lambda t_n}\sigma^\top(t_n,x_n)q^u_n)|
\end{align*}
Repeating the above argument and using that $f$ is jointly continuous in $(t,x)$ uniformly in $(y,z)$ we get that the upper limit of the right-hand side when $n\to\infty$ is bounded by $k_f(u(t_n,x_n)-w(t_n,y_n))$. Put together, this gives that
\begin{align*}
(\lambda-k_f) \limsup_{n\to\infty}(u(t_n,x_n)-w(t_n,y_n))&\leq 0
\end{align*}
and choosing $\lambda>k_f$ gives a contradition.\qed\\ 

\bibliographystyle{plain}
\bibliography{qvi-rbsde_ref}
\end{document}